\documentclass{article}

\usepackage{geometry}
\usepackage{fancyhdr}
\usepackage{amsmath, amssymb}
\usepackage{multirow}     
\usepackage[toc,page]{a  ppendix}
\usepackage{placeins}
\usepackage{graphicx}
\usepackage{soul}
\usepackage{sans}
\usepackage{array}
\newcommand{\PreserveBackslash}[1]{\let\temp=\\#1\let\\=\temp}
\newtheorem{thm}{Theorem}

\newtheorem{Conj}[thm]{Conjecture}
\newcolumntype{C}[1]{>{\PreserveBackslash\centering}p{#1}}
\newcolumntype{R}[1]{>{\PreserveBackslash\raggedleft}p{#1}}
\newcolumntype{L}[1]{>{\PreserveBackslash\raggedright}p{#1}}
\usepackage{subfig}
\usepackage[numbers]{natbib}
\usepackage{titling}
\usepackage{textcomp}
\usepackage{color}
\definecolor{bg}{rgb}{0.89, 0.95, 0.71} 
\definecolor{ac}{rgb}{0.50, 0.18, 0.41} 
\definecolor{rc}{rgb}{0.77, 0.57, 0.71} 
\usepackage{listings}
\newcommand{\ts}{\textsuperscript}
\usepackage{wrapfig}
\usepackage{subfig}

\RequirePackage{lineno}
\modulolinenumbers[1]

\usepackage{natbib}
\makeatletter
\newcommand{\specificthanks}[1]{\@fnsymbol{#1}}
\makeatother

\title{Accuracy Enhancing Interface Treatment Algorithm: the Back and Forth Error Compensation and Correction method} 

\author{ 
Wenbin Dong\thanks{wdong000@citymail.cuny.edu, Civil Engineering Department, City College of New York, City University of New York,
New York, NY 10031, USA.} 
\and 
Yingjie Liu\thanks{yingjie@math.gatech.edu,
School of Mathematics, Georgia Institute of Technology, Atlanta, GA
30332.} 
\and 
Hansong Tang\thanks{htang@ccny.cuny.edu, Civil Engineering Department, City College of New York, City University of New York,
New York, NY 10031, USA.}
}

\begin{document}

\maketitle

\renewcommand{\thefootnote}{\fnsymbol{footnote}}

\footnotetext{Declarations:
{\bf Funding}. This work is supported in part by National Science Foundation (DMS $\#$ 1622459 and 1622453). 
{\bf Conflicts of interest/Competing interests}. There is no conflicts of interest.
{\bf Availability of data and material}. Data is available within the article.
{\bf  Code availability}. The code is available upon request.
}

\begin{abstract}
The accuracy of information transmission while solving domain decomposed problems is crucial to smooth transition of a solution around the interface/overlapping region. This paper describes a systematical study on an accuracy enhancing interface treatment algorithm based on the back and forth error compensation and correction method (BFECC). By repetitively employing a low order interpolation technique (normally local 2\ts{nd} order) 3 times, this algorithm achieves local 3\ts{rd} order accuracy. Analytical derivations for 1D \& 2D cases are made, and the ``super convergence'' phenomenon (4\ts{th} order accuracy) is found for specific positioning of the donor and target grids. A set of numerical experiments based on various relative displacements, relative rotations, mesh ratios, and meshes with perturbation have been tested, and the results match the derivations. Different interface treatments are compared with 3D examples: corner flow and cavity flow. The component convergence rate analysis shows that the BFECC method positively affects the accuracy of solutions.
\end{abstract}

\section{Introduction}

\vspace{0.2cm} 

It is important to have accurate numerical treatments at grid or model interfaces and for those of the algorithms to be coupled in subdomains. Commonly the algorithms within each subdomain are 2\ts{nd} order accurate, or, locally 3\ts{rd} order accurate, while interface algorithms are at most locally 2\ts{nd} order accurate, e.g., in flow simulations  \cite{Benek1983,Tang_2014_b,Tang1999}. Therefore, it is important to develop locally 3\ts{rd} order accurate methods which are easy to use on nonuniform grids for interface algorithms, so that the numerical accuracy will be  uniformly 2\ts{nd} order in the whole computational domain. \par

Interpolation is a key step for interface treatments, and there exist many high-order interpolation methods, such as Lagrangian methods and spline methods, e.g., \cite{Chesshire1990,Epperson_2013_b}. However, in general, a high-order interpolation needs 1) more grid nodes whose values are interpolated and 2) solving a linear system to determine interpolation coefficients (in particular for irregular grids \cite{Barth_1990_b}). When the spatial dimension gets higher, both 1) and 2) are expensive for computation, and they are often difficult to implement in practical computation, e.g., in interpolation between irregular grids. To achieve 3\ts{rd} order accuracy locally at interfaces requires looking beyond adjacent grid points by existing methods. \par

A high order conservative interpolation method for moving meshes was introduced in \cite{Tang_2003_b} which proposed a fictitious advection equation representing the mesh movement and solved it with a high order finite volume scheme to achieve the conservative interpolation. In this research, we utilize the BFECC method \cite{Liuyj_2003_b} to achieve 3\ts{rd} order accuracy locally by calling a local linear interpolation algorithm $3$ times.  For example, for an advection equation ($u$ is a scalar, and ${\bf V}$ is a vector)
\begin{equation}
\label{advec_eq}
u_t +{\bf V}\cdot \bigtriangledown u = 0~,
\end{equation}

\noindent with the 1\ts{st} order CIR scheme \cite{Courant_1952_b} as the underlying scheme, BFECC provides a procedure, as outlined in Sec.2, to enhance the accuracy to 2\ts{nd} order accurate, thus improving the local order of accuracy to 3\ts{rd} order. The advantage of the BFECC method is that it does not increase the complexity of the simple underlying scheme, but just repeat the computation using the same scheme for two extra times. \par

When there is no advection, or, ${\bf V}=0$ in Eq. (\ref {advec_eq}), the grid at an old-time level is the original grid and the grid at a new-time level is the grid to be interpolated to. Then,  the BFECC method becomes an interpolation method that achieves local 3\ts{rd} order accuracy if a local linear interpolation is used in the underlying scheme. In comparison with standard higher-order accurate interpolation methods, BFECC keeps the simplicity (and three times the complexity) of the underlying linear interpolation, and it can be implemented to existing computer codes by recycling an existing 2\ts{nd} order interface algorithms. There has been research work on BFECC for quadtree-octree mesh refinement,  e.g., in \cite{kim2005flowfixer,kim2007advections} for fluid simulations, \cite{TNOV15} for transport equations, \cite{Yu_2011_b} for locally refined Lattice Boltzmann method,  and for interpolation between equilateral triangular grids \cite{Hu_2014_b}. 
The interpolation can also be viewed as solving Eq. (\ref {advec_eq}) on the fixed old grid by a numerical scheme with ${\bf V}$ describing the mesh movement between the old and new grids \cite{Tang_2003_b}, and 
the case with the numerical scheme being BFECC was mentioned in \cite{Liuyj_2007_b}. 

\par

In this paper,  the BFECC method is analyzed and experimented with regard to its implementation, effectiveness, and limitations as a systematical study on its application to interpolation on structured meshes. For a comparison, 2\ts{nd} order interpolation and MFBI (mass-flux based interpolation) are also discussed \cite{Tang2003}. The rest of the paper is organized as follows: in Section~\ref{BFECC_int}, we introduce the algorithm of BFECC and discuss its accuracy as well as some other concerns in 1D, 2D \& 3D cases. Section~\ref{Sec: flow} provides the numerical simulation results for cavity flow and corner flow. Then we finish this paper in Section~\ref{Sec: Conclusion} with conclusions.

\section{BFECC Interpolation}
\label{BFECC_int}

\subsection{Algorithm and Accuracy}
\label{subsec: Alg}
\vspace{0.2cm} 

Considering equation (\ref{advec_eq}) in $d$ space dimensions as an example, Let $L$
be a numerical scheme for (\ref{advec_eq}) which advances the numerical
solution at the time level $t_n$, $U^n$ forward to the 
next time level $t_{n+1}$, {\it i.e.}, $U^{n+1}=L\{ U^n\} $.
The BFECC method \cite{Liuyj_2003_b} can be described as
follows. \par

\begin{enumerate} 
\item {\bf Forward advection}:

$ \widetilde{U}^{n+1} = L\{U^n\}~.$

\item {\bf Backward advection}:

$ \widetilde{U}^n = L^*\{\widetilde{U}^{n+1}\}~,$

where $L^*$ is $L$ applied to the time reversed equation of (\ref{advec_eq}), $u_t -{\bf V}\cdot \bigtriangledown u = 0$.

\item {\bf Error compensation}:

$ \hat{U}^n = U^n+\dfrac{U^n-\widetilde{U}^n}{2}~.$

\item {\bf Forward advection}:

$ U^{n+1} = L\{ \hat{U}^n \}~.$
\end{enumerate}

Assume scheme $L$ is a linear scheme on a rectangular grid $\mathcal{G}_1 $ with mesh sizes $\Delta x_1= \Delta x_2= \cdots =\Delta x_d\equiv \Delta x$, and the Fourier symbols of $L$ and $L^*$, $\rho_L$ and $\rho_{L^*}$ respectively,
satisfy  $\rho_{L^*}=\overline{\rho_L}$ (the complex conjugate of $\rho_L$), then we have the
following theorem (\cite{Liuyj_2007_b}, Theorem~4).
\begin{thm}
\label{thm: acc}
Suppose scheme $L$ is
accurate of order $r$ for equation~(\ref{advec_eq}) with constant coefficients, 
where $r$ is an odd positive integer, then after applying BFECC to
scheme $L$, it is accurate of order $r+1$.
\end{thm}

Let $\mathcal{G}_2 $ be another rectangular grid generated by
shifting grid $\mathcal{G}_1 $ along an arbitrary constant vector ${\bf w}$, $|{\bf w}|=O(\Delta x)$. Let $f$ be a sufficiently smooth function defined
at grid points of grid $\mathcal{G}_1 $, then interpolating
the value of $f$ from grid points of $\mathcal{G}_1 $ to grid points of $\mathcal{G}_2 $ can be performed as follows analogous to the above BFECC method. 

\begin{enumerate} 

\item {\bf Forward interpolation}: 

Denote the value of $f$ at grid point ${\bf i}$ as $f_{\bf i}$, for every grid point ${\bf i}$ of $\mathcal{G}_1 $. Use local linear interpolation (e.g., bilinear interpolation in 2D, trilinear interpolation in 3D, etc) to interpolate the value of $f$ from grid points of $\mathcal{G}_1 $ to every grid point of $\mathcal{G}_2 $.

\item {\bf Backward interpolation}:

Use local linear interpolation to interpolate from grid points of $\mathcal{G}_2 $ to grid point ${\bf i}$ of $\mathcal{G}_1 $ to obtain $\tilde{f_{\bf i}}$, for every grid point ${\bf i}$ of $\mathcal{G}_1 $.

\item {\bf Error compensation}:

Define $\hat{f_{\bf i}}=f_{\bf i} + \frac12 (f_{\bf i}-\tilde{f_{\bf i}})$, for every grid point ${\bf i}$ of $\mathcal{G}_1 $.

\item {\bf Forward interpolation}:

Use local linear interpolation to interpolate the values $\{\hat{f_{\bf i}}\}$ from grid points of $\mathcal{G}_1 $ to every grid point of $\mathcal{G}_2 $.

\end{enumerate}

This interpolation method will be referred to as {\bf BFECC
interpolation} for the rest of the paper. Note that the above BFECC interpolation can be applied to interpolation between two irregular or unstructured meshes as long as the local linear interpolation is properly defined, e.g., as local linear least squares.
Now we are going to
identify BFECC interpolation with solving an advection equation
with BFECC following \cite{Tang_2003_b}.
Replace ${\bf V}$
with ${\bf -w}/\Delta x $ in equation (\ref{advec_eq}) and let the 
initial value of this equation be $f$. Then solve equation (\ref{advec_eq}) for just one time step with the time step size $\Delta t=\Delta x$ by 
the scheme formed by BFECC applied to the CIR scheme \cite{Courant_1952_b} (with local linear interpolation). It is easy to see that
the numerical solution at grid point ${\bf i}$ and at the time $\Delta t$ is the interpolated value at grid point ${\bf i+w}$
of grid $\mathcal{G}_2 $ obtained by the above BFECC
interpolation. BFECC applied to the CIR scheme is 2\ts{nd} order accurate according to
Theorem~\ref{thm: acc}, or equivalently is 3\ts{rd} order accurate for just one time step. We have proved the following theorem.

\begin{thm}
\label{thm:BFECC_interp}
BFECC interpolation from a uniform rectangular grid in $k$ dimensions (for any positive integer $k$, with mesh sizes $\Delta x_1=\Delta x_2=\cdots=\Delta x_k$) to a shifted grid 
along any vector is 3\ts{rd} order accurate.
\end{thm}

Since the interpolation between two grids is equivalent to solving an advection equation for one time step by a numerical scheme, we have also tested another simple scheme, the modified MacCormack scheme \cite{Selle_2008_b} which uses only two advection steps. Note that Theorem~\ref{thm:BFECC_interp} can be generalized when the local linear interpolation is replaced by
local interpolation with a $(2k+1)$-th degree polynomial, $k=1,2,3,\cdots$, and the corresponding order of BFECC interpolation becomes $2k+3$.






\subsection{Super-Convergence} 

\begin{figure}[!htbp]
\centering
\includegraphics[scale=0.6]{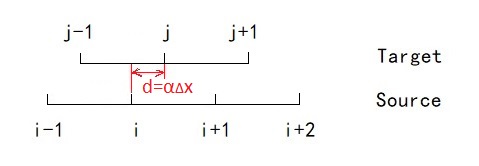}
\caption{1D grids for interpolation}
\label{fig:1d_prove}
\end{figure}

Consider in one dimension the BFECC interpolation for a sufficiently smooth function $f(x)$, from $f_i$ defined on the source grid to $f_j$ on the target grid, as shown in Fig.\ref{fig:1d_prove}. Here the target grid is obtained
by shifting the source grid (uniform with grid size $\Delta x$) to the right with a distance $d$. Let $\alpha = d/\Delta x$.
The function value at node $j$ of the target grid is interpolated from its values in the source grid as follows:

\noindent Forward interpolation:

\begin{equation}
  \begin{array}{ll}
  \cdots \\
  {f}^*_{j-1}= \alpha f_i   + (1-\alpha)f_{i-1} \\
  {f}^*_j    = \alpha f_{i+1} + (1-\alpha)f_i   \\
  {f}^*_{j+1}= \alpha f_{i+2} + (1-\alpha)f_{i+1} \\
  \cdots
  \end{array}
\end{equation}

\noindent Backward interpolation: 

\begin{equation}
  \begin{array}{ll}
   \cdots \\
  \tilde{f}_i= \alpha f^*_{j-1}   + (1-\alpha)f^*_j \\
  \tilde{f}_{i+1} = \alpha f^*_j + (1-\alpha)f^*_{j+1} \\
   \cdots 
  \end{array}
\end{equation}

\noindent Error compensation: 

\begin{equation}
  \begin{array}{ll}
   \cdots \\
\hat{f}_i  = f_i   + (f_i-\tilde{f}_i)/2 \\
\hat{f}_{i+1}= f_{i+1} + (f_{i+1}-\tilde{f}_{i+1})/2 \\
\cdots
 \end{array}
\end{equation}

\noindent Forward interpolation:
\begin{equation}
  \begin{array}{ll}
  \cdots \\ 
f^{new}_j = \alpha \hat{f}_{i+1} + (1-\alpha)\hat{f}_i   \\
 \cdots 
  \end{array}
\end{equation}
Here $f^{new}_j $ is the interpolated value at node $j$ 
by BFECC interpolation.
Using Taylor expansions (at node $j$), a straightforward derivation shows that
\begin{equation}
\label{1d_BFECC_Taylor}
f^{new}_j = f_j + \frac{\alpha(\alpha-1)[(8\alpha-4)\Delta x^3(f_{xxx})_j -9\alpha(\alpha-1)\Delta x^4(f_{xxxx})_j]}{24}+{\cal O}(\Delta x^4)~.
\end{equation}

Formula~(\ref{1d_BFECC_Taylor}) is derived using a program of symbolic calculation in the software of 'Mathematica' (version: 11.3.0). The code is attached in Appendix \ref{1D_BFECC_proof}. This implies that the error of the BFECC interpolation is 3\ts{rd} order accurate. {\bf And if $\alpha= \frac12$, {\it i.e.}, the interpolation point is located at the centroid of a cell of the grid, the term with $\Delta x^3$ on the right hand side disappears, then the interpolation has $4$th order accuracy, a ``super-convergence'' phenomenon.} As a comparison, if only the forward interpolation step is applied, then

\begin{equation}
\label{1d_int_Taylor}
f^*_j = f_j + \frac{\alpha(1-\alpha)\Delta x^2 (f_{xx})_j}{2}+{\cal O}(\Delta x^3)~,
\end{equation}

\noindent which is 2\ts{nd} order accurate. 

Now consider a 2D rectangular grid $\mathcal{G}_1=\{(i\Delta x, j\Delta y): \; {\rm for}\;{\rm any}\;{\rm integers}\; i,j\}$. Let $\mathcal{G}_2$ be the grid generated by shifting $\mathcal{G}_1$ along the vector $(\alpha \Delta x, \beta \Delta y) $, 
$0\le \alpha, \beta\le 1$,
see Fig.\ref {fig:2d_prove}. Consider the BFECC interpolation from $\mathcal{G}_1$ to $\mathcal{G}_2$ for a sufficiently smooth function $f(x,y)$.  Using 2D Taylor expansions around an interpolation point $(r,s)=(i+\alpha, j+\beta)$, we obtain the interpolation accuracy similar to the 1D case (\ref {1d_BFECC_Taylor}):   

\begin{figure}[!htbp]
\centering
\includegraphics[scale=0.5]{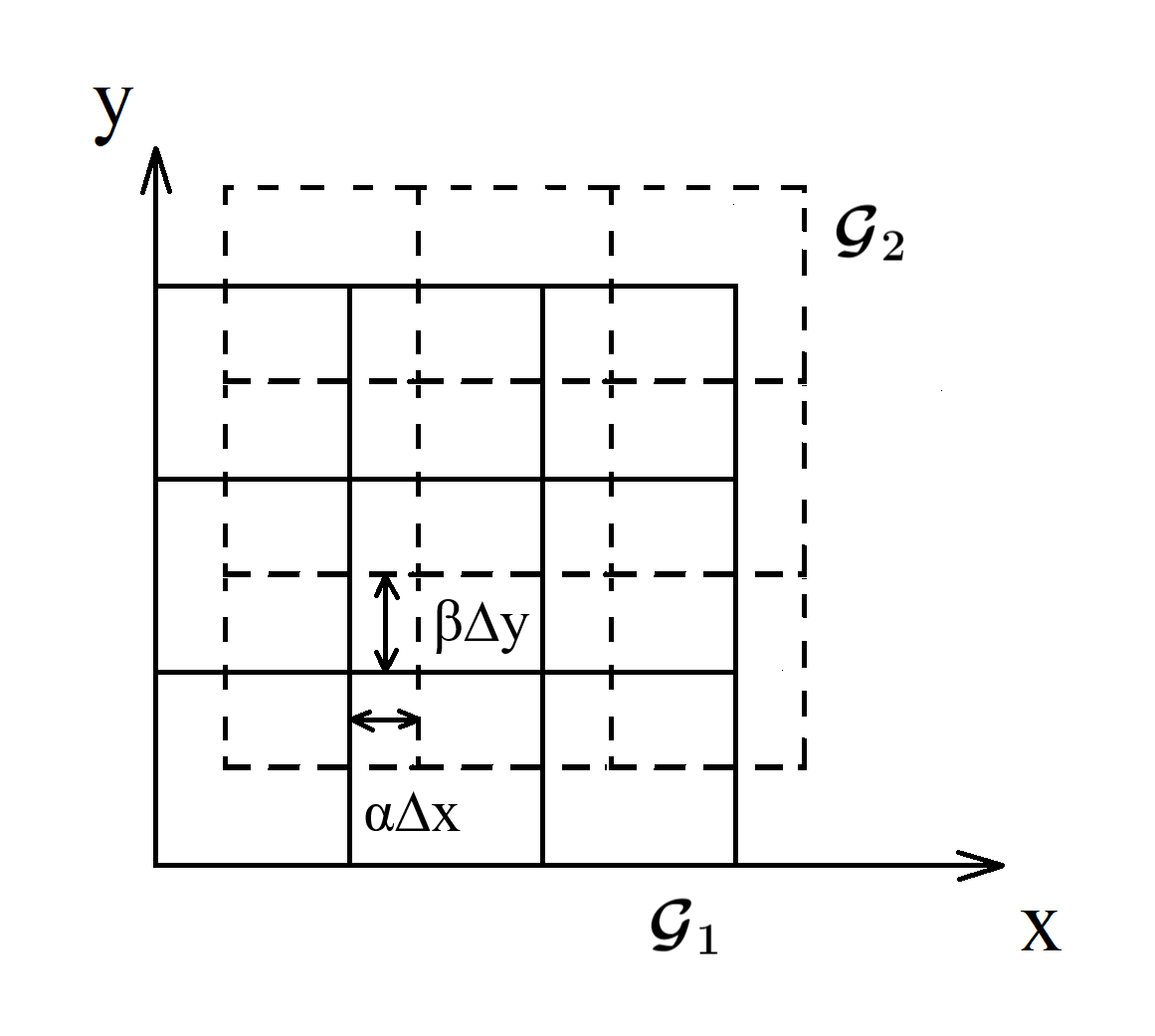}
\caption{2D grids for interpolation}
\label{fig:2d_prove}
\end{figure}

\begin{equation}
\label{2d_BFECC_Taylor}
\begin{array}{ll}
f_{(r,s)}^{new} = f_{(r,s)}&+\dfrac{(\alpha-3\alpha^2+2\alpha^3)\Delta x^3(f_{xxx})_{(r,s)}+(\beta-3\beta^2+2\beta^3)\Delta y^3(f_{yyy})_{(r,s)}}{6} \\
         &+{\cal O}(\Delta x^4+\Delta y^4),
\end{array}
\end{equation}
   
\noindent where $f_{(r,s)}^{new}$ is the interpolated value at
shifted grid point $(r,s)$ which is 3\ts{rd} order accurate. 
{\bf In particular, when $\alpha=\beta=\frac12$, $f_{(r,s)}^{new}$ is
4\ts{th} order accurate.} The 'Mathematica' (version: 11.3.0) code to derive Formula~(\ref{2d_BFECC_Taylor})  is attached in Appendix \ref{2D_BFECC_proof}. As a comparison, bi-linear interpolation at $(r,s)$ has

\begin{equation}
\label{2d_int_Taylor}
\begin{array}{ll}
f_{(r,s)}^*= f_{(r,s)} &+ \frac12 (\alpha-\alpha^2)\Delta x^2(f_{xx})_{(r,s)} \\
       &+ \frac12 (\beta-\beta^2)\Delta y^2(f_{yy})_{(r,s)}+{\cal O}(\Delta x^3+\Delta y^3)~,
\end{array}
\end{equation}

which is only 2\ts{nd} order accurate. We summarize these results
in the following theorem.

\begin{thm}
\label{thm:BFECC_interp_2D}
BFECC interpolation from a uniform rectangular grid in 1D or
2D (with mesh sizes $\Delta x$ and $\Delta y$) to a shifted grid 
along any vector is 3\ts{rd} order accurate. In particular, if the
shifted grid points are located at centroids of the original grid cells, then BFECC interpolation is 4\ts{th} order accurate.
\end{thm}

Numerical experiments show that Theorem~\ref{thm:BFECC_interp_2D} holds in 3D. We have the following conjecture.

\begin{Conj}
\label{thm:BFECC_interp_d_dim}
BFECC interpolation from a uniform rectangular grid in $d$ dimensions (for any positive integer $d$, with mesh sizes $\Delta x_1, \Delta x_2, \cdots, \Delta x_d$) to a shifted grid 
along any vector is 3\ts{rd} order accurate. In particular, if the
shifted grid points are located at centroids of the original grid cells, then BFECC interpolation is 4\ts{th} order accurate.
\end{Conj}

\subsection{Numerical Tests} 

In order to verify the above discussion, numerical tests have been conducted.  To obtain the order of convergence rate using linear interpolation and BFECC interpolation (based on the linear interpolation.), four sets of grids are used with
grid sizes $\Delta x=0.05$, $0.025$, $0.0125$ \& $0.00625$ . The computed orders of accuracy against a 1D exact function are listed in Table~\ref{tab:1d_linear_parabolic_bfecc_ca}.  The 2D convergence results can be found in Table~\ref{tab:2d_bilinear_bfecc_a}.
When interpolation points are located at centroids of grid cells,
BFECC interpolation gives 4\ts{th} order accuracy in 2D, see Table~\ref{tab:2d_bilinear_bfecc_a_super}. In Table~\ref{tab:2d_lsq_bfecc_mac_1}, BFECC interpolation is compared to the interpolation based on the modified MacCormack scheme (the interpolation is performed by solving an advection equation with
the scheme for one time step, see the discussion near the end of Sec.~\ref{subsec: Alg}.) While both methods have 3\ts{rd} order accuracy, BFECC interpolation errors tend to be several times smaller. Similar comparison tests are conducted with different interpolated functions in Tables~\ref{tab:2d_lsq_bfecc_mac_2} and
\ref{tab:2d_lsq_bfecc_mac_3}, and with interpolation points located at centroids of grid cells in Table~\ref{tab:2d_lsq_bfecc_mac_4}. Note that BFECC interpolation has 4\ts{th} order accuracy at cenroids while the interpolation based on the modified MacCormack scheme doesn't have ``super-convergence'' there. 
Bilinear interpolation is replaced by linear least squares method in some examples as the underlying interpolation method which doesn't change the order of convergence. Least squares interpolation is suitable for irregular meshes, see e.g., \cite{Barth_1990_b, hu1999weighted} for high order least squares interpolation for WENO reconstruction on unstructured meshes and \cite{wangBackForthError2019} for using linear least squares on irregular meshes to construct an underlying scheme of
BFECC to solve the Maxwell's equations. We will extend BFECC interpolation to non uniform grids next and linear least squares interpolation will be used as the underlying interpolation method for BFECC. When the interpolation point is at the centroid of a rectangle, both bilinear interpolation and linear least squares are simple average of function values at vertices of the rectangle. 

\begin {table}[!htbp]
\begin{center}
\begin{tabular}{| C{2cm} | C{2cm} | C{2cm} | C{2cm} | C{2cm} |}
\hline
\multicolumn{5}{|c|}{1D, Grid Shift ${\bf w}=0.25\Delta x$} \\ \hline
  & \multicolumn{2}{c}{Linear Interpolation} & \multicolumn{2}{|c|}{BFECC} \\ \hline
$\Delta$x & error   & order    & error   & order    \\ \hline
0.05    & 7.19E-4 &      & 5.80E-5 &      \\ \hline
0.025   & 1.67E-4 & 2.11 & 7.26E-6 & 3.00 \\ \hline
0.0125  & 4.00E-5 & 2.06 & 9.10E-7 & 3.00 \\ \hline
0.00625 & 9.81E-6 & 2.03 & 1.14E-7 & 3.00 \\ \hline
\end{tabular}
\caption {Accuracy of linear interpolation and BFECC interpolation in 1D. Interpolated function: $f(x)=sin(\pi*x)$.}
\label{tab:1d_linear_parabolic_bfecc_ca}
\end{center}
\end {table}

\begin {table}[!htbp]
\begin{center}
\begin{tabular}{| C{2cm} | C{2cm} | C{2cm} | C{2cm} | C{2cm} |}
\hline
\multicolumn{5}{|c|}{2D, Grid Shift ${\bf w}=(0.25\Delta x,0.25\Delta y)$, $\Delta x=\Delta y$} \\ \hline
  & \multicolumn{2}{c}{Bilinear Interpolation} & \multicolumn{2}{|c|}{BFECC} \\ \hline
$\Delta$x & error   & order    & error   & order    \\ \hline
0.05    & 2.26E-3 &      & 8.40E-5 &      \\ \hline
0.025   & 5.24E-4 & 2.11 & 1.08E-5 & 2.95 \\ \hline
0.0125  & 1.27E-4 & 2.05 & 1.33E-6 & 3.03 \\ \hline
0.00625 & 3.11E-5 & 2.03 & 1.18E-7 & 2.99 \\ \hline
\end{tabular}
\caption {Accuracy of bilinear and BFECC interpolations in 2D. Interpolated function: $f(x,y)=sin(\pi*(x+2y))$.}
\label{tab:2d_bilinear_bfecc_a}
\end{center}
\end {table}
  
\begin {table}[!htbp]
\begin{center}
\begin{tabular}{| C{2cm} | C{2cm} | C{2cm} | C{2cm} | C{2cm} |}
\hline
\multicolumn{5}{|c|}{2D, Grid Shift ${\bf w}=(0.5\Delta x,0.5\Delta y)$, $\Delta x=\Delta y$} \\ \hline
  & \multicolumn{2}{c}{Bilinear Interpolation} & \multicolumn{2}{|c|}{BFECC} \\ \hline
$\Delta$x  & error & order & error & order \\ \hline
0.05    & 1.93E-3 &      & 2.01E-4 &      \\ \hline
0.025   & 4.57E-4 & 2.08 & 1.27E-5 & 3.98 \\ \hline
0.0125  & 1.12E-4 & 2.03 & 7.93E-7 & 4.01 \\ \hline
0.00625 & 2.76E-5 & 2.01 & 4.98E-8 & 3.99 \\ \hline
\end{tabular}
\caption {Super-convergence of BFECC interpolation based on bilinear interpolation in 2D. Interpolated function: $f(x,y)=sin(\pi*(x+2y))$.}
\label{tab:2d_bilinear_bfecc_a_super}
\end{center}
\end {table}

\begin {table}[!htbp]
\begin{center}
\begin{tabular}{| C{1.5cm} | C{1.5cm} | C{1.5cm} | C{1.5cm} | C{1.5cm} | C{1.5cm} | C{1.5cm} |}
\hline
\multicolumn{7}{|c|}{2D, Grid Shift ${\bf w}=(0.25\Delta x,0.25\Delta y)$, $\Delta x=\Delta y$} \\ \hline
  & \multicolumn{2}{c}{Linear least squares} & \multicolumn{2}{|c|}{BFECC} & \multicolumn{2}{|c|}{Maccormack} \\ \hline
$\Delta$x & error & order & error & order & error & order \\ \hline
0.05    & 7.63E-4 &      & 5.85E-5  &      & 1.66E-4  &      \\ \hline
0.025   & 1.76E-4 & 2.11 & 7.31E-6  & 3.00 & 2.08E-5  & 2.99 \\ \hline
0.0125  & 4.24E-5 & 2.06 & 9.16E-7  & 3.00 & 2.61E-6  & 3.00 \\ \hline
0.00625 & 1.04E-5 & 2.03 & 1.15E-7  & 3.00 & 3.27E-7  & 3.00 \\ \hline
\end{tabular}
\caption {Accuracy of linear least squares, and BFECC and Maccormack interpolation methods based on it in 2D. Interpolated function: $f(x,y)=sin(\pi*(x+0.2y))$.}
\label{tab:2d_lsq_bfecc_mac_1}
\end{center}
\end {table}

\begin {table}[!htbp]
\begin{center}
\begin{tabular}{| C{1.5cm} | C{1.5cm} | C{1.5cm} | C{1.5cm} | C{1.5cm} | C{1.5cm} | C{1.5cm} |}
\hline
\multicolumn{7}{|c|}{2D, Grid Shift ${\bf w}=(0.25\Delta x,0.25\Delta y)$, $\Delta x=\Delta y$} \\ \hline
  & \multicolumn{2}{c}{Linear least squares} & \multicolumn{2}{|c|}{BFECC} & \multicolumn{2}{|c|}{Maccormack} \\ \hline
$\Delta$x & error & order & error & order & error & order \\ \hline
0.05    & 1.53E-4 &      & 1.18E-5  &      & 3.38E-5  &      \\ \hline
0.025   & 3.51E-5 & 2.12 & 1.48E-6  & 3.00 & 4.22E-6  & 3.00 \\ \hline
0.0125  & 8.42E-6 & 2.06 & 1.85E-7  & 3.00 & 5.27E-7  & 3.00 \\ \hline
0.00625 & 2.06E-6 & 2.03 & 2.31E-8  & 3.00 & 6.59E-8  & 3.00 \\ \hline
\end{tabular}
\caption {Accuracy of linear least squares, and BFECC and Maccormack interpolation methods based on it in 2D. Interpolated function: $f(x,y)=(x+0.2y)^3$.}
\label{tab:2d_lsq_bfecc_mac_2}
\end{center}
\end {table}

\begin {table}[!htbp]
\begin{center}
\begin{tabular}{| C{1.5cm} | C{1.5cm} | C{1.5cm} | C{1.5cm} | C{1.5cm} | C{1.5cm} | C{1.5cm} |}
\hline
\multicolumn{7}{|c|}{2D, Grid Shift ${\bf w}=(0.25\Delta x,0.25\Delta y)$, $\Delta x=\Delta y$} \\ \hline
  & \multicolumn{2}{c}{Linear least squares} & \multicolumn{2}{|c|}{BFECC} & \multicolumn{2}{|c|}{Maccormack} \\ \hline
$\Delta$x & error & order & error & order & error & order \\ \hline
0.05    & 2.52E-4 &      & 1.93E-6  &      & 5.71E-6  &      \\ \hline
0.025   & 6.22E-5 & 2.02 & 2.45E-7  & 2.98 & 7.11E-7  & 3.01 \\ \hline
0.0125  & 1.55E-5 & 2.01 & 3.08E-8  & 2.99 & 8.87E-8  & 3.00 \\ \hline
0.00625 & 3.85E-6 & 2.00 & 3.86E-9  & 2.99 & 1.11E-8  & 3.00 \\ \hline
\end{tabular}
\caption {Accuracy of linear least squares, and BFECC and Maccormack interpolation methods based on it in 2D. Interpolated function: $f(x,y)=exp(x+0.2y)$.}
\label{tab:2d_lsq_bfecc_mac_3}
\end{center}
\end {table}

\begin {table}[!htbp]
\begin{center}
\begin{tabular}{| C{1.5cm} | C{1.5cm} | C{1.5cm} | C{1.5cm} | C{1.5cm} | C{1.5cm} | C{1.5cm} |} \hline
\multicolumn{7}{|c|}{2D, Grid Shift ${\bf w}=(0.5\Delta x,0.5\Delta y)$, $\Delta x=\Delta y$} \\ \hline
  & \multicolumn{2}{c}{Linear least squares} & \multicolumn{2}{|c|}{BFECC} & \multicolumn{2}{|c|}{Maccormack} \\ \hline
$\Delta$x=$\Delta$y  & error & order & error & order & error & order \\ \hline
0.05    & 3.38E-4 &      & 1.65E-7  &      & 9.95E-6  &      \\ \hline
0.025   & 8.32E-5 & 2.02 & 1.01E-8  & 4.02 & 1.24E-6  & 3.01 \\ \hline
0.0125  & 2.06E-5 & 2.01 & 6.29E-10 & 4.01 & 1.54E-7  & 3.00 \\ \hline
0.00625 & 5.14E-6 & 2.01 & 3.91E-11 & 4.01 & 1.92E-8  & 3.00 \\ \hline
\end{tabular}
\caption {Super-convergence study of linear least squares, and BFECC and Maccormack interpolation methods based on it in 2D. Interpolated function: $f(x,y)=exp(x+0.2y)$.}
\label{tab:2d_lsq_bfecc_mac_4}
\end{center}
\end {table}

\newpage

BFECC interpolation on 3D rectangular grids is also tested in Tables~\ref{tab:3d_trilinear_bfecc_ca_025} and \ref{tab:3d_lsq_bfecc_a},
with the underlying interpolation method being trilinear and linear least squares respectively. Again, ``super-convergence'' is
found for BFECC interpolation if interpolation points are located at centroids of grid cells, see Table~\ref{tab:3d_trilinear_bfecc_ca_050}.
We have also tested 3D uniform rectangular grids with $\Delta x\ne \Delta y\ne\Delta z $ and found Theorem~\ref{thm:BFECC_interp_2D} still true.
These numerical tests are presented in Tables~\ref{tab:3d_tri_nonuniform_025}, \ref{tab:3d_tri_nonuniform_050} and \ref{tab:3d_lsq_nonuniform_025}. 

\begin {table}[!htbp]
\begin{center}
\begin{tabular}{| C{2cm} | C{2cm} | C{2cm} | C{2cm} | C{2cm} |} \hline
\multicolumn{5}{|c|}{3D, Grid Shift ${\bf w}=(0.25\Delta x,0.25\Delta y, 0.25\Delta z)$, $\Delta x=\Delta y=\Delta z$} \\ \hline
  & \multicolumn{2}{c}{Trilinear Interpolation} & \multicolumn{2}{|c|}{BFECC} \\ \hline
$\Delta$x & error & order    & error      & order    \\ \hline
0.05    & 2.28E-2 &      & 1.88E-3 &      \\ \hline
0.025   & 5.68E-3 & 2.00 & 2.06E-4 & 3.19 \\ \hline
0.0125  & 1.41E-3 & 2.01 & 2.48E-5 & 3.05 \\ \hline
0.00625 & 3.51E-4 & 2.01 & 3.08E-6 & 3.01 \\ \hline
\end{tabular}
\caption {Accuracy of trilinear and BFECC interpolation in 3D. Interpolated function: $f(x,y,z)=sin(\pi*(x+2y+3z))$.}
\label{tab:3d_trilinear_bfecc_ca_025}
\end{center}
\end {table}

\begin {table}[!htbp]
\begin{center}
\begin{tabular}{| C{2cm} | C{2cm} | C{2cm} | C{2cm} | C{2cm} |}
\hline
\multicolumn{5}{|c|}{3D, Grid Shift ${\bf w}=(0.25\Delta x,0.25\Delta y, 0.25\Delta z)$, $\Delta x=\Delta y=\Delta z$} \\ \hline
  & \multicolumn{2}{c}{Linear Least Squares} & \multicolumn{2}{|c|}{BFECC} \\ \hline
$\Delta$x & error & order    & error   & order    \\ \hline
0.05    & 1.12E-2 &      & 1.03E-3 &      \\ \hline
0.025   & 2.73E-3 & 2.04 & 1.29E-4 & 2.99 \\ \hline
0.0125  & 6.73E-3 & 2.02 & 1.64E-5 & 2.99 \\ \hline
0.00625 & 1.67E-4 & 2.01 & 2.04E-6 & 2.99 \\ \hline
\end{tabular}
\caption {Accuracy of linear least squares and BFECC interpolation based on it in 3D. Interpolated function: $f(x,y,z)=sin(\pi*(x+2y+3z))$.}
\label{tab:3d_lsq_bfecc_a}
\end{center}
\end {table}

\begin {table}[!htbp]
\begin{center}
\begin{tabular}{| C{2cm} | C{2cm} | C{2cm} | C{2cm} | C{2cm} |}
\hline
\multicolumn{5}{|c|}{3D, Grid Shift ${\bf w}=(0.5\Delta x,0.5\Delta y, 0.5\Delta z)$, $\Delta x=\Delta y=\Delta z$} \\ \hline
  & \multicolumn{2}{c}{Trilinear Interpolation} & \multicolumn{2}{|c|}{BFECC} \\ \hline
$\Delta$x & error & order    & error   & order    \\ \hline
0.05    & 3.02E-2 &      & 1.90E-3 &      \\ \hline
0.025   & 7.57E-3 & 2.13 & 1.22E-4 & 3.97 \\ \hline
0.0125  & 1.88E-3 & 2.06 & 7.60E-6 & 4.00 \\ \hline
0.00625 & 4.68E-4 & 2.03 & 4.74E-7 & 4.00 \\ \hline
\end{tabular}
\caption {Super-convergence of BFECC interpolation based on trilinear interpolation in 3D. Interpolated function: $f(x,y,z)=sin(\pi*(x+2y+3z))$.}
\label{tab:3d_trilinear_bfecc_ca_050}
\end{center}
\end {table}

\begin {table}[!htbp]
\begin{center}
\begin{tabular}{| C{2cm} | C{2cm} | C{2cm} | C{2cm} | C{2cm} |} \hline
\multicolumn{5}{|c|}{3D, Grid Shift ${\bf w}=(0.25\Delta x,0.25\Delta y, 0.25\Delta z)$, $\Delta x:\Delta y:\Delta z$=1:0.9:1.2} \\ \hline
  & \multicolumn{2}{c}{Trilinear Interpolation} & \multicolumn{2}{|c|}{BFECC} \\ \hline
$\Delta$x & error & order    & error      & order    \\ \hline
0.05    & 2.85E-2 &      & 2.75E-3 &      \\ \hline
0.025   & 7.21E-3 & 1.99 & 2.97E-4 & 3.21 \\ \hline
0.0125  & 1.80E-3 & 2.00 & 3.55E-5 & 3.06 \\ \hline
0.00625 & 4.49E-4 & 2.00 & 4.39E-6 & 3.01 \\ \hline
\end{tabular}
\caption {Accuracy of trilinear and BFECC interpolation in 3D. Interpolated function: $f(x,y,z)=sin(\pi*(x+2y+3z))$.}
\label{tab:3d_tri_nonuniform_025}
\end{center}
\end {table}

\begin {table}[!htbp]
\begin{center}
\begin{tabular}{| C{2cm} | C{2cm} | C{2cm} | C{2cm} | C{2cm} |} \hline
\multicolumn{5}{|c|}{3D, Grid Shift ${\bf w}=(0.5\Delta x,0.5\Delta y, 0.5\Delta z)$, $\Delta x:\Delta y:\Delta z$=1:0.9:1.2} \\ \hline
  & \multicolumn{2}{c}{Trilinear Interpolation} & \multicolumn{2}{|c|}{BFECC} \\ \hline
$\Delta$x & error & order    & error      & order    \\ \hline
0.05    & 3.77E-2 &      & 2.90E-3 &      \\ \hline
0.025   & 9.59E-3 & 1.97 & 1.89E-4 & 3.94 \\ \hline
0.0125  & 2.40E-3 & 2.00 & 1.19E-5 & 3.99 \\ \hline
0.00625 & 5.99E-4 & 2.00 & 7.44E-7 & 4.00 \\ \hline
\end{tabular}
\caption {Super-convergence of BFECC interpolation based on trilinear interpolation in 3D. Interpolated function: $f(x,y,z)=sin(\pi*(x+2y+3z))$.}
\label{tab:3d_tri_nonuniform_050}
\end{center}
\end {table}

\begin {table}[!htbp]
\begin{center}
\begin{tabular}{| C{2cm} | C{2cm} | C{2cm} | C{2cm} | C{2cm} |} \hline
\multicolumn{5}{|c|}{3D, Grid Shift ${\bf w}=(0.25\Delta x,0.25\Delta y, 0.25\Delta z)$, $\Delta x:\Delta y:\Delta z$=1:0.9:1.2} \\ \hline
  & \multicolumn{2}{c}{Linear Least Squares} & \multicolumn{2}{|c|}{BFECC} \\ \hline
$\Delta$x & error & order    & error      & order    \\ \hline
0.05    & 1.59E-2 &      & 8.26E-4 &      \\ \hline
0.025   & 3.92E-3 & 2.02 & 8.86E-5 & 3.22 \\ \hline
0.0125  & 9.73E-4 & 2.01 & 1.06E-5 & 3.06 \\ \hline
0.00625 & 2.42E-4 & 2.01 & 1.32E-6 & 3.01 \\ \hline
\end{tabular}
\caption {Accuracy of linear least squares and BFECC interpolation in 3D. Interpolated function: $f(x,y,z)=sin(\pi*(x+2y+3z))$.}
\label{tab:3d_lsq_nonuniform_025}
\end{center}
\end {table}

\newpage

\subsection{Interpolation Between Grids of Different Spacing}  

\begin{figure}[!htbp]
 \centering
  \subfloat{\label{ratio1}\includegraphics[width=0.35\textwidth]{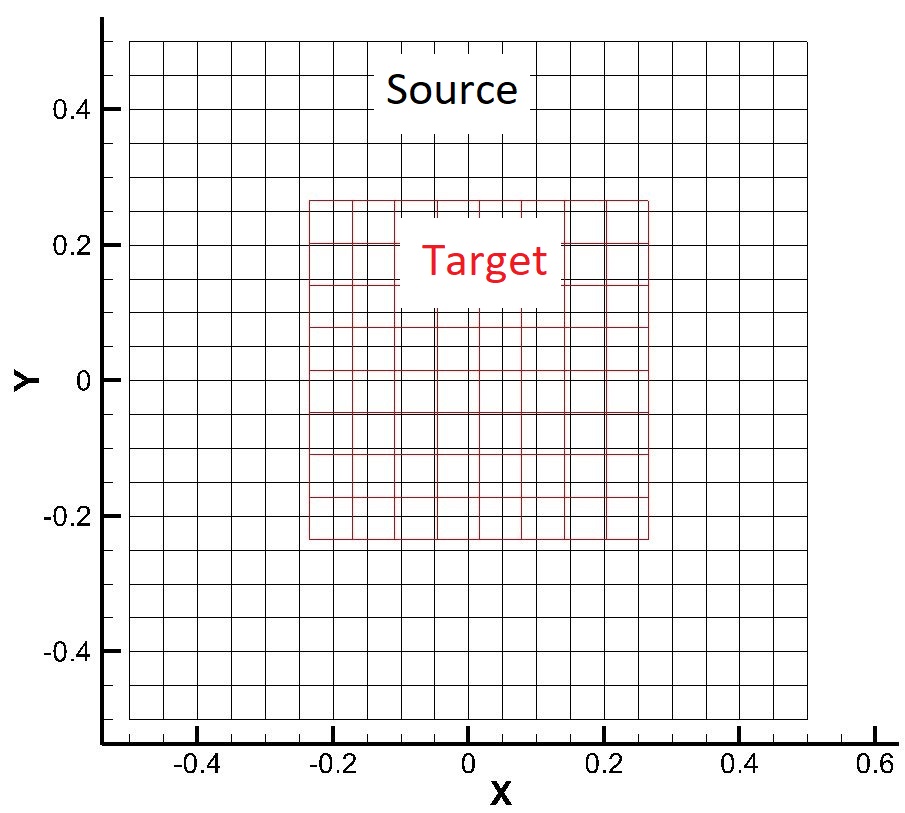}}  \ \
  \subfloat{\label{ratio2}\includegraphics[width=0.35\textwidth]{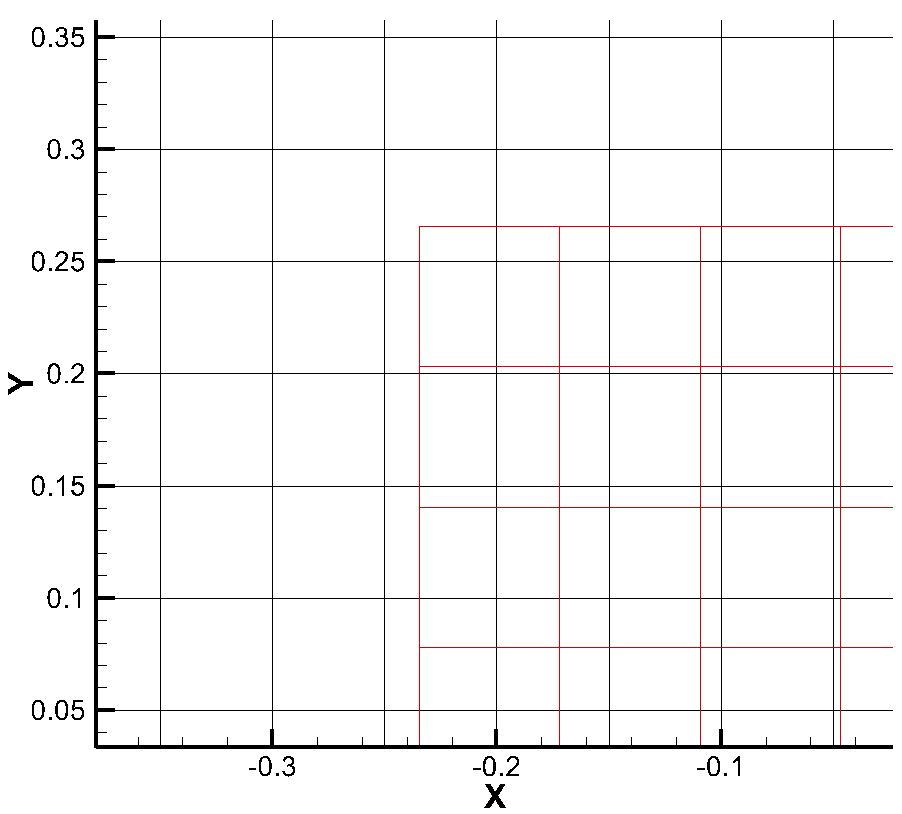}}
  \caption{Grid size ratio of $\sqrt{2}:1$.}
  \label{fig:ratio_mesh}
\end{figure}

BFECC interpolation can be implemented between two non-uniform grids with ease. We consider the interpolation accuracy when the BFECC interpolation
is between two parallel rectangular grids of different spacing. A 3D numerical experiment is conducted as shown in Fig.\ref {fig:ratio_mesh}, in which BFECC interpolation 
(with the underlying trilinear interpolation) will be made from grid $S$ to $T$ with grid sizes $\sqrt{2}:1$ respectively. 
The results are presented in Table~\ref {tab:3d_tri_bfecc_a_ratio}. It is clear that, in comparison to trilinear interpolation, BFECC interpolation reduces its errors in general but makes no improvement of the convergence order. 

\begin {table}[!htbp]
\begin{center}
\begin{tabular}{| C{2cm} | C{2cm} | C{2cm} | C{2cm} | C{2cm} |}
\hline
\multicolumn{5}{|c|}{3D, Grid spacing ratio of $\sqrt{2}:1$} \\ \hline
  & \multicolumn{2}{c}{Trilinear} & \multicolumn{2}{|c|}{BFECC} \\ \hline
$\Delta$x & error    & order    & error   & order    \\ \hline
0.0625    & 5.21E-3 &      & 2.89E-3 &      \\ \hline
0.03125   & 1.21E-3 & 2.11 & 6.95E-4 & 2.05 \\ \hline
0.015625  & 2.67E-4 & 2.18 & 1.74E-4 & 2.00 \\ \hline
0.0078125 & 5.96E-5 & 2.16 & 4.24E-5 & 2.04 \\ \hline
\end{tabular}
\caption {Accuracy \& order of trilinear and BFECC interpolation for 3D grids of size ratio $\sqrt{2}:1$. Interpolated function: $f(x,y,z)=sin(\pi*(x+2y+3z))$.}
\label{tab:3d_tri_bfecc_a_ratio}
\end{center}
\vspace{-10pt}
\end {table}

\begin{figure}[!htbp]
\centering
\includegraphics[scale=0.2]{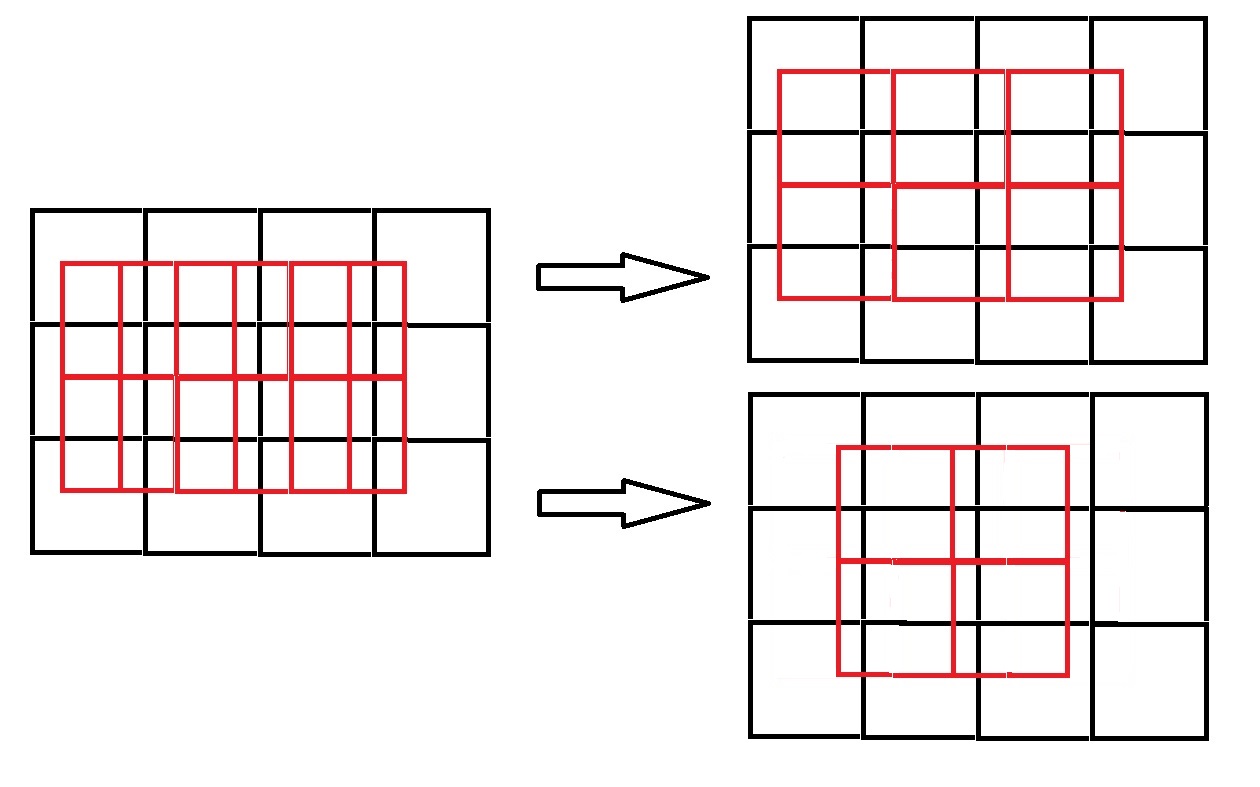}
\caption{Fine grid is divided into 2 coarse grids to meet the 1:1 ratio.}
\label{fig:ratio2}
\end{figure}

If the BFECC interpolation is from a coarse grid to a fine grid and the ratio of two grid spacing is far from 1:1, a remedy is to view the fine grid as the overlap of several subgrids so the ratio of the grid size of the coarse grid to that of a subgrid is close to 1:1. Fig.~\ref {fig:ratio2} shows an example which has 1:2 grid size ratio of a fine grid, the target grid, and a coarse grid, the source grid, and the splitting of the fine grid into two subgrids. Then the BFECC interpolation can be applied to the interpolation between the coarse grid (the source grid) and each of the subgrids to achieve 3\ts{rd} order accuracy. An numerical example of the 3D BFECC interpolation between the coarse grid and a subgrid in the fine grid similar to Fig.\ref{fig:ratio2}. 

\subsection{Smoothly Perturbed Rectangular Grids}  
\begin{figure}[!htbp]
 \centering
  \subfloat{\label{perturbed1}\includegraphics[width=0.3\textwidth]{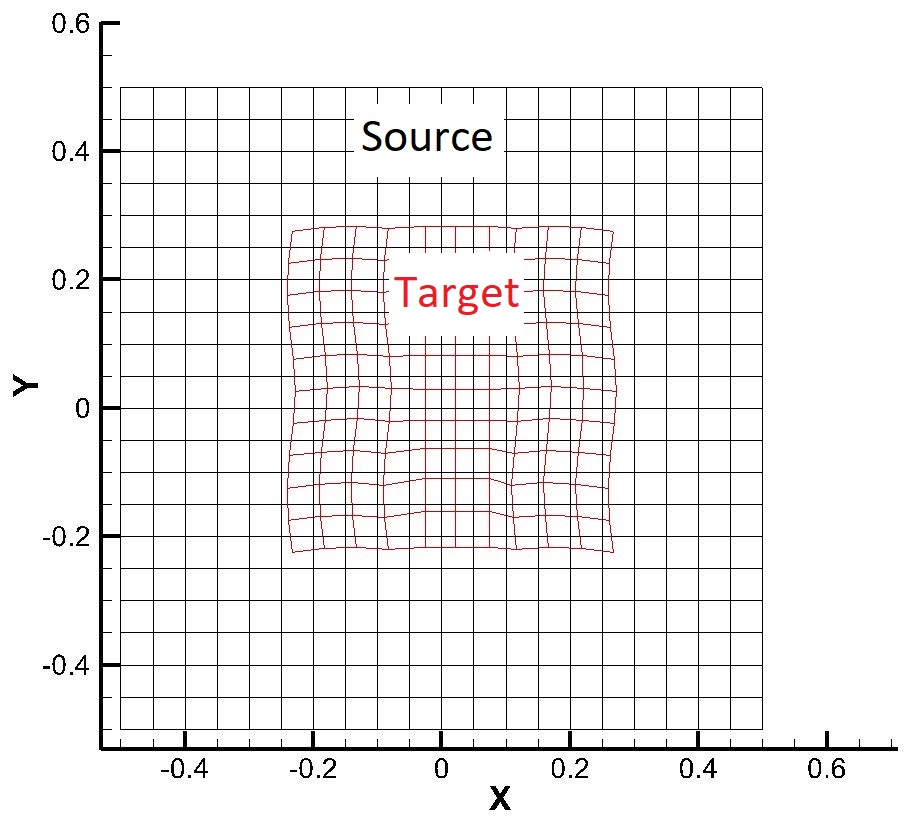}}  \ \
  \subfloat{\label{perturbed2}\includegraphics[width=0.3\textwidth]{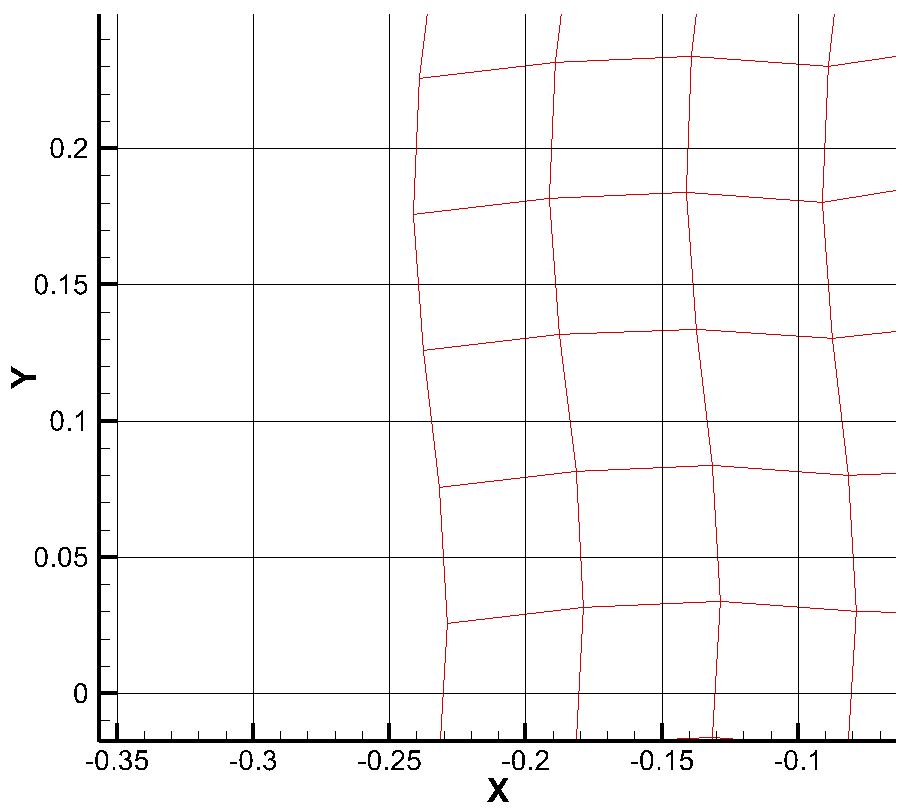}}
  \caption{The target grid is smoothly perturbed and shifted from the source grid.}
  \label{fig:perturbed_mesh}
\end{figure}

 BFECC interpolation from a rectangular grid to a shifted and smoothly perturbed grid still has 3\ts{rd} order accuracy. For example,  consider the grid defined by (\ref{un_mesh}) which is the shift and smooth perturbation from a rectangular grid ((\ref{un_mesh}) without the terms containing ``$\sin$'' functions), see Fig.~\ref{fig:perturbed_mesh}.
 
\begin{equation}
\label{un_mesh}
\begin{array}{ll}
x(i,j,k)=(i-1)\Delta x - 0.5L_x + 0.1\Delta x\; sin(\pi*((i-1)\Delta x - 0.5L_x)), \\
y(i,j,k)=(j-1)\Delta y - 0.5L_y + 0.1\Delta y\; sin(\pi*((j-1)\Delta y - 0.5L_y)), \\
z(i,j,k)=(k-1)\Delta z - 0.5L_z + 0.1\Delta z\; sin(\pi*((k-1)\Delta z - 0.5L_z)),
\end{array}
\end{equation}

\noindent where '$L_x$', '$L_y$', and '$L_z$'  are lengths of the domain in x, y and z directions. The target grid (the former one) deviates up to 20\% from the the source grid (the latter one).  The BFECC interpolation of the function $ f=sin(x+2y+z^2) $ from the source grid ($\Delta x=\Delta y= \Delta z$) to the target grid is shown in Table~\ref {tab:3d_trilinear_bfecc_ca_perturbed}. With the non uniform deformation from the source grid, the order of accuracy of BFECC interpolation still gets around 3\ts{rd} order.  

\begin {table}[!htbp]
\begin{center}
\begin{tabular}{| C{2cm} | C{2cm} | C{2cm} | C{2cm} | C{2cm} |}
\hline
\multicolumn{5}{|c|}{From rectangular grid with $\Delta x=\Delta y= \Delta z$ to smoothly perturbed grid} \\ \hline
\multicolumn{5}{|c|}{Grid Shift $w=(0.25\Delta x,0.25\Delta y,0.25\Delta z)$} \\ \hline
  & \multicolumn{2}{c}{Trilinear} & \multicolumn{2}{|c|}{BFECC} \\ \hline
$\Delta$x & error   & order    & error   & order    \\ \hline
0.0125    & 6.26E-5 &      & 8,81E-6 &      \\ \hline
0.00625   & 1.42E-5 & 2.14 & 1.27E-6 & 2.79 \\ \hline
0.003125  & 3.31E-6 & 2.10 & 1.70E-7 & 2.91 \\ \hline
0.0015625 & 7.97E-7 & 2.05 & 2.04E-8 & 3.06 \\ \hline
\end{tabular}
\caption {Accuracy of BFECC interpolation on perturbed and shifted grid in 3D. }
\label{tab:3d_trilinear_bfecc_ca_perturbed}
\end{center}
\end {table}

\subsection{Rotated Rectangular Grids}

When the grid is rotated, the BFECC interpolation from the original grid to it tends to improve the accuracy of the underlying interpolation method, but no longer improves the order of convergence. Here is a test. After the rotation of a source grid with an angle, the results of BFECC interpolation are presented in Table~\ref{tab:3d_trilinear_bfecc_ca_rotated}. This test shows that the rotation has prevented the accuracy order of the BFECC interpolation from reaching 3\ts{rd} order, but still the errors of BFECC interpolations is several times smaller than those of the underlying trilinear interpolation.

\begin {table}[!htbp]
\begin{center}
\begin{tabular}{| C{2cm} | C{2cm} | C{2cm} | C{2cm} | C{2cm} |}
\hline
\multicolumn{5}{|c|}{Rotation of $3^{\circ}$ of a Rectangular grid} \\ \hline
  & \multicolumn{2}{c}{Trilinear} & \multicolumn{2}{|c|}{BFECC} \\ \hline
$\Delta$x & error & order    & error   & order    \\ \hline
0.05    & 3.45E-3 &      & 1.51E-3 &      \\ \hline
0.025   & 7.48E-4 & 2.21 & 1.62E-4 & 3.22 \\ \hline
0.0125  & 2.03E-4 & 1.88 & 3.56E-5 & 2.18 \\ \hline
0.00625 & 4.61E-5 & 2.14 & 1.26E-5 & 1.50 \\ \hline
\end{tabular}
\caption {Accuracy of BFECC interpolation on rotated grid in 3D.}
\label{tab:3d_trilinear_bfecc_ca_rotated}
\end{center}
\end {table}


\section{Application to Flow Problems} 
\label{Sec: flow}
\subsection{Corner Flow}  

A 3D corner flow was proposed to test accuracy of numerical schemes \cite{Tang_2007_b}, and it is selected to test the performance of BFECC in simulation of flow problems. Consider a cubic box with length of 1.2 in each direction. Let $u$, $v$, and $w$ be velocity in $x$, $y$, and $z$, respectively, and the origin of the coordinates are at the center of the cube). The fluid inside is initialized with zero speed, then an entrance velocity is imposed at the lid of the box as below:
\begin{equation}
\label{entrance_vel}
w=-e^{1-2/t}sin(\pi((x+0.5)^3-3(x+0.5)^2+3(x+0.5))(0.25-y^2))
\end{equation}

\noindent It gradually drives the fluid to flow out at one of the lateral sides. The Reynolds number is Re=100, based on the edge length and the maximum velocity at the lid. The flow is simulated by solving the incompressible Navier-Stokes equations on overset grids together with trilinear interpolation at grid interfaces \cite{Tang_2006_b}. The Navier-Stokes equations are discretized by a 2\ts{nd} order finite difference method in both time and space.  

\begin{figure}[!htbp]
\centering
\includegraphics[scale=0.3]{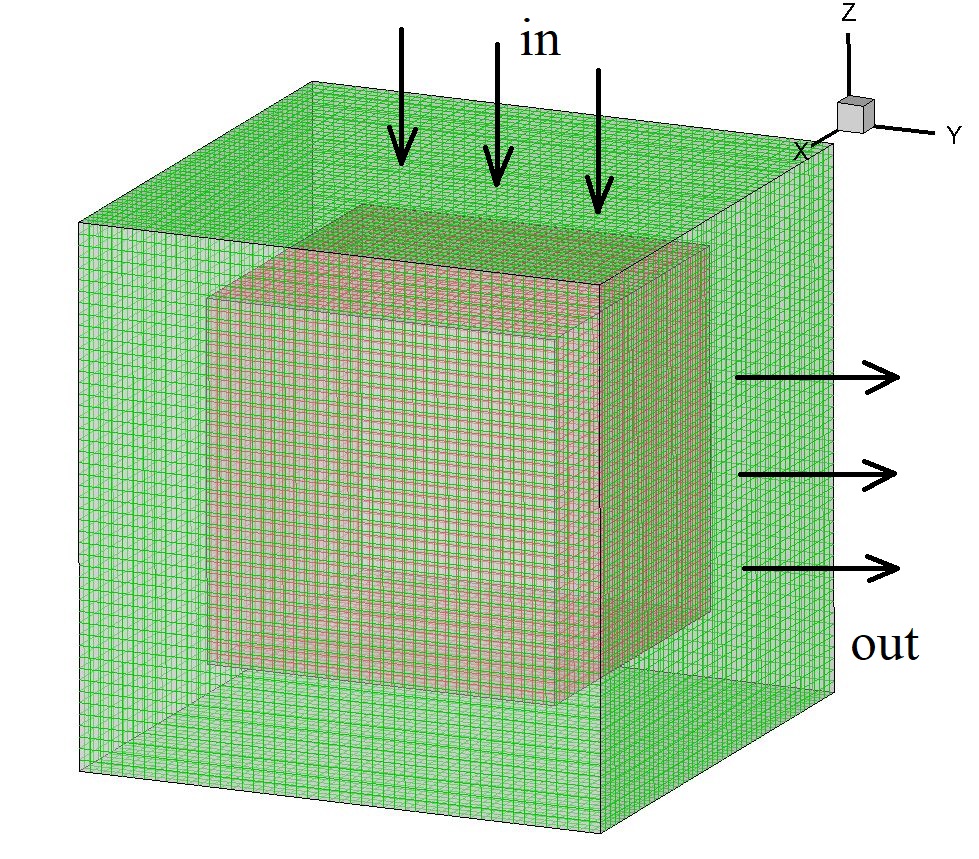}
\caption{\small Fine mesh for the corner flow}
\label{fig:corner_mesh}
\end{figure}

There are two overlapping cubes, the small cube (0.8*0.8*0.8) is placed inside the big cube (1.2*1.2*1.2). The edges of two cubes are parallel but the nodes are not coincided. Three sets of grids are used: $13^3$ for the big cube and $9^3$ for the small one as coarse grids, $25^3$ and $17^3$ as medium grids, and $49^3$ and $33^3$ as fine grids, see Fig.\ref{fig:corner_mesh}. Additionally, a simulation is made used a single domain of $71^3$, and its solution is intended as a reference to be compared with.  

\begin{figure}[!htbp]
 \centering
  \subfloat[pressure (y=0)]{\label{corner_hor_p}\includegraphics[width=0.35\textwidth]{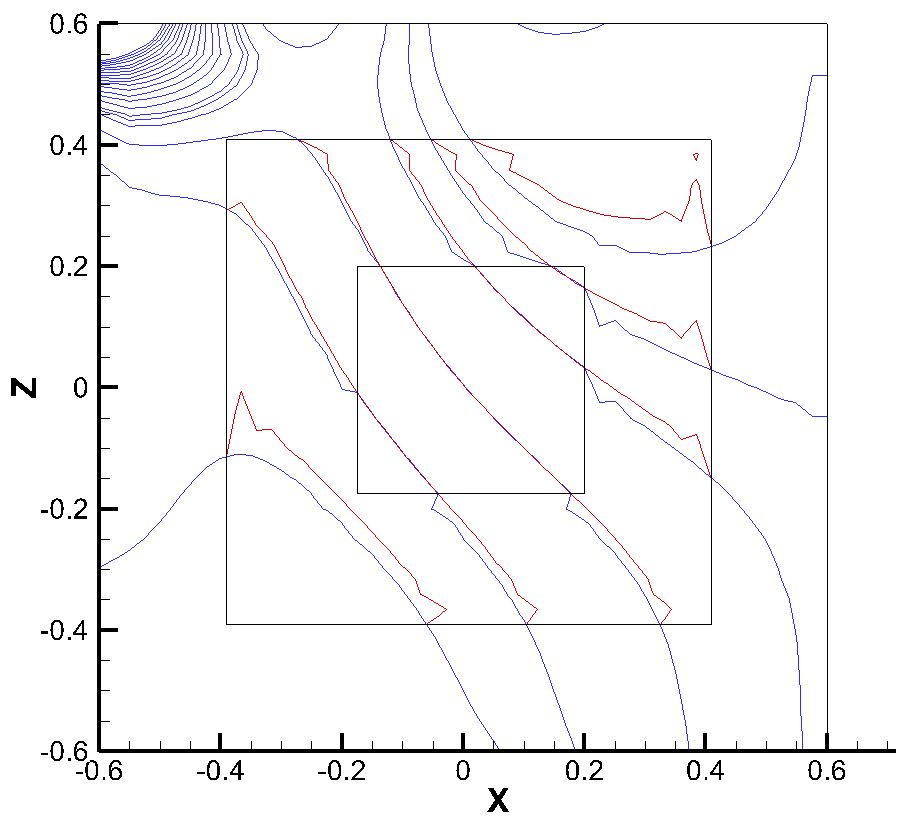}}  \ \
  \subfloat[u (y=0)]{\label{corner_hor_u}\includegraphics[width=0.35\textwidth]{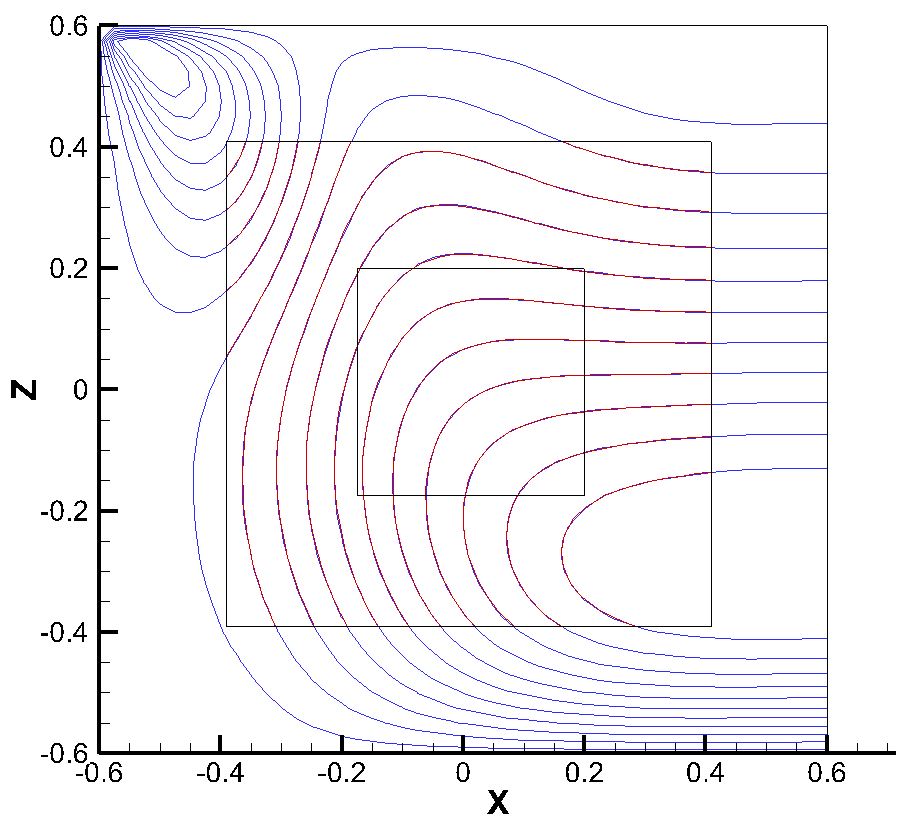}}\\
  \subfloat[v (x=0)]{\label{corner_hor_v}\includegraphics[width=0.35\textwidth]{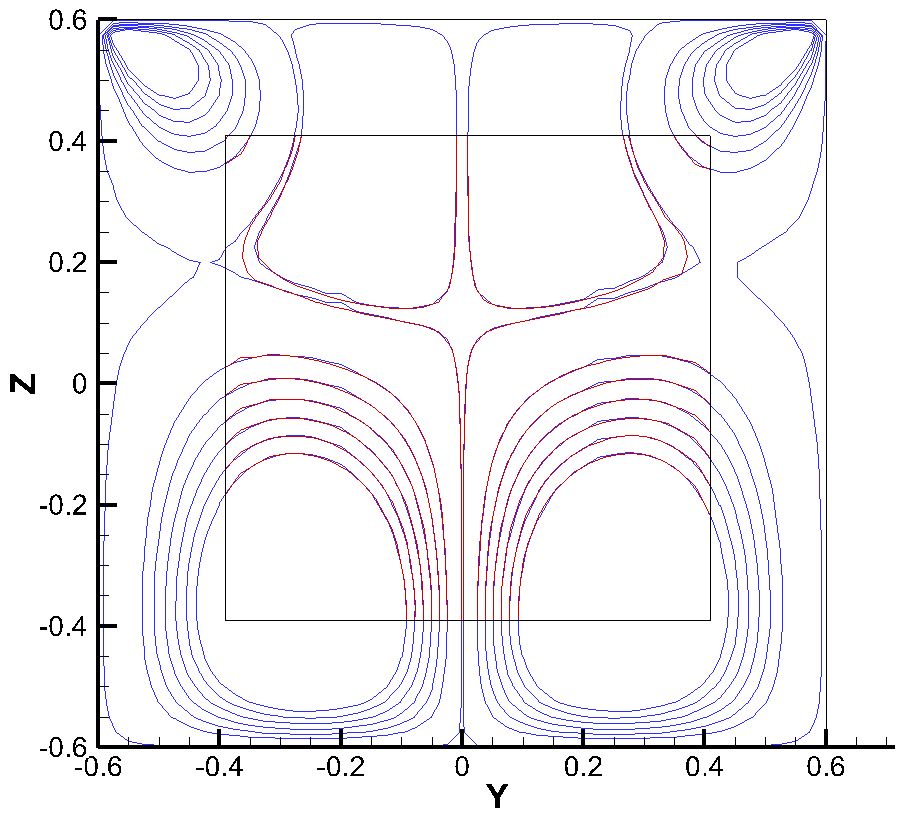}}  \ \
  \subfloat[w (y=0)]{\label{corner_hor_w}\includegraphics[width=0.35\textwidth]{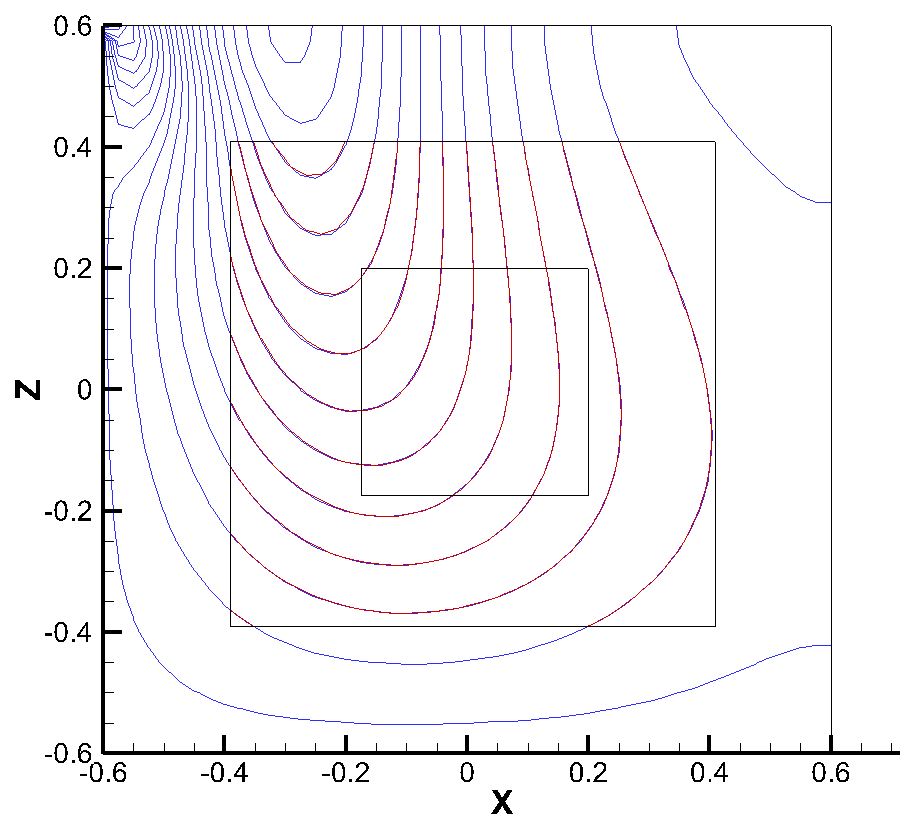}}\\  
  \caption{\small Contour line at t=10s. The interfaces are treated by BFECC.}
  \label{fig:corner_hor for BFECC}
\end{figure}

Simulated instantaneous flow fields are shown are shown in Fig. \ref {fig:corner_hor for BFECC}. It is seen that the solutions on the outer and inner cubes are identical at the scales of the figures. Further check on the solutions are made by looking at their accuracy orders with regard to time step and grid spacing. As stated in the Lax Theorem, in linear situations, an algorithm of $pth$-order leads to $pth$-order in accuracy order. In non-linear situation, such understanding is frequently true. In view that the algorithms to solve the flows inside cubes are 2\ts{nd} order accurate, and thus their convergence rates are expected to be 2\ts{nd} order. However, when the interface computation is treated by linear interpolation, which is 1\ts{st} order, the convergence rate will be around 1\ts{st} order. As a result, the overall accuracy order of solutions to the corner flows will be affected, or, degraded by the linear interpolation. As BFECC is applied at the interfaces, it is expected that the accuracy order will be recovered to 2\ts{nd} order. 

With above understanding and also in view that there is no analytical solution for the flow, the accuracy order, $\kappa $, of the solutions are evaluated numerically by their values on three sets of grids listed above using the following formula (e.g., \cite{Tang_2007_b}): 

\begin{equation}
\label{order_estimate}
\kappa =log(\dfrac{||f_{\Delta }-f_{\Delta /2}||_2}{||f_{\Delta /2}-f_{\Delta /4}||_2})/log2
\end{equation}

\noindent where $f$ is solution for flow variables $u$, $v$, $w$, and $p$, and $\Delta  $, $\Delta  /2$, and $\Delta  /4$ stand for coarse, medium, and fine grids, respectively. It is noted that formula (\ref {order_estimate}) is a rough estimate for the accuracy order and may not present its exact value. 

\begin{figure}[!htbp]
 \centering
  \subfloat[p]{\label{hor_inter_ca_p}\includegraphics[width=0.35\textwidth]{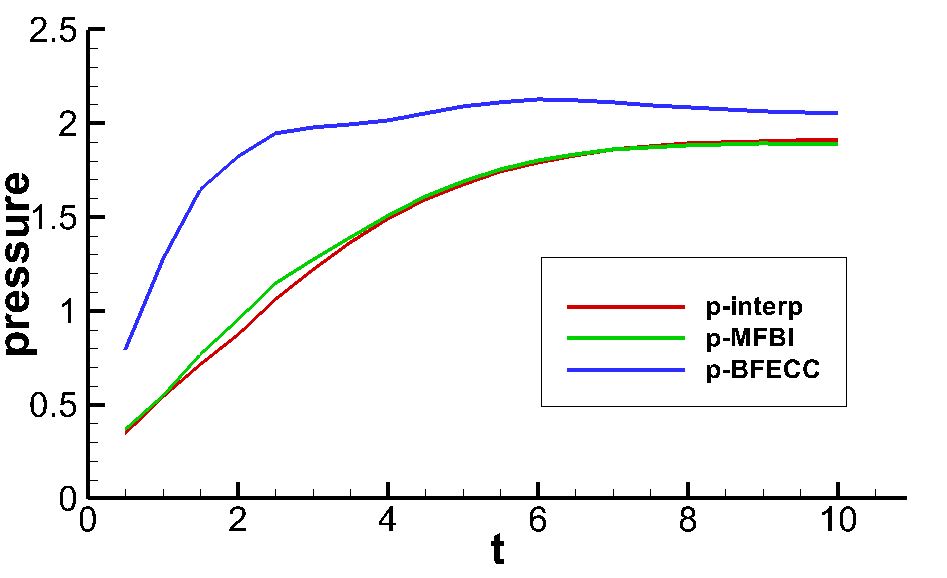}}  \ \
  \subfloat[u]{\label{hor_inter_ca_u}\includegraphics[width=0.35\textwidth]{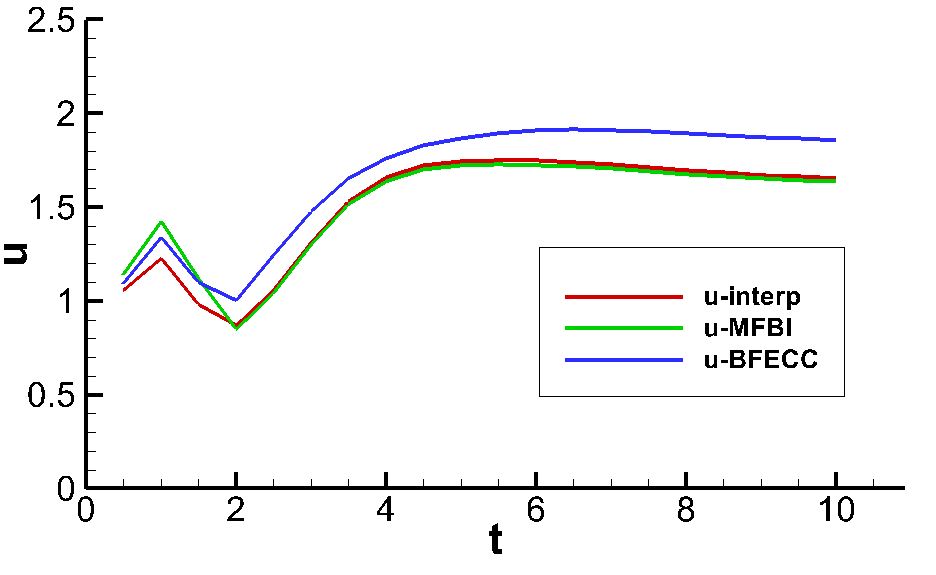}} \\
  \subfloat[v]{\label{hor_inter_ca_v}\includegraphics[width=0.35\textwidth]{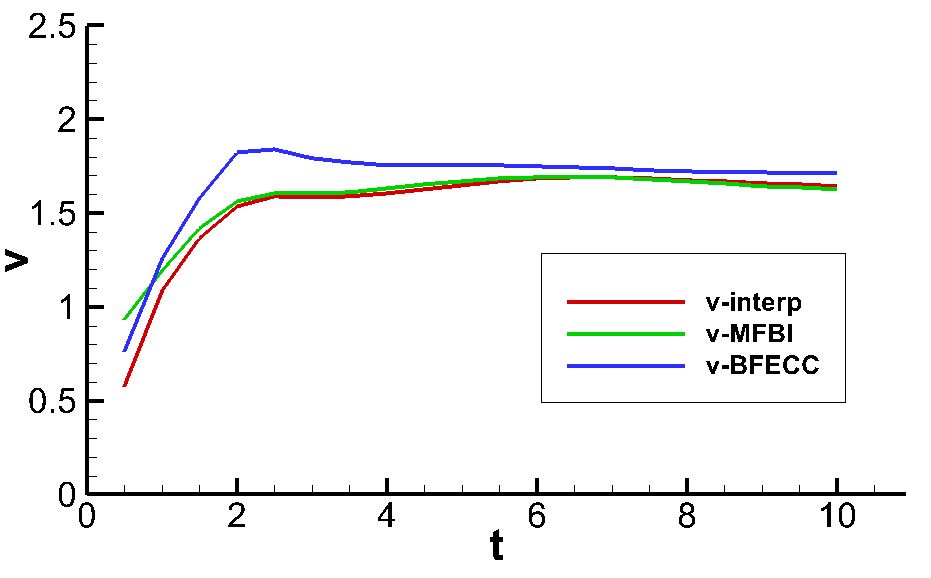}}  \ \
  \subfloat[w]{\label{hor_inter_ca_w}\includegraphics[width=0.35\textwidth]{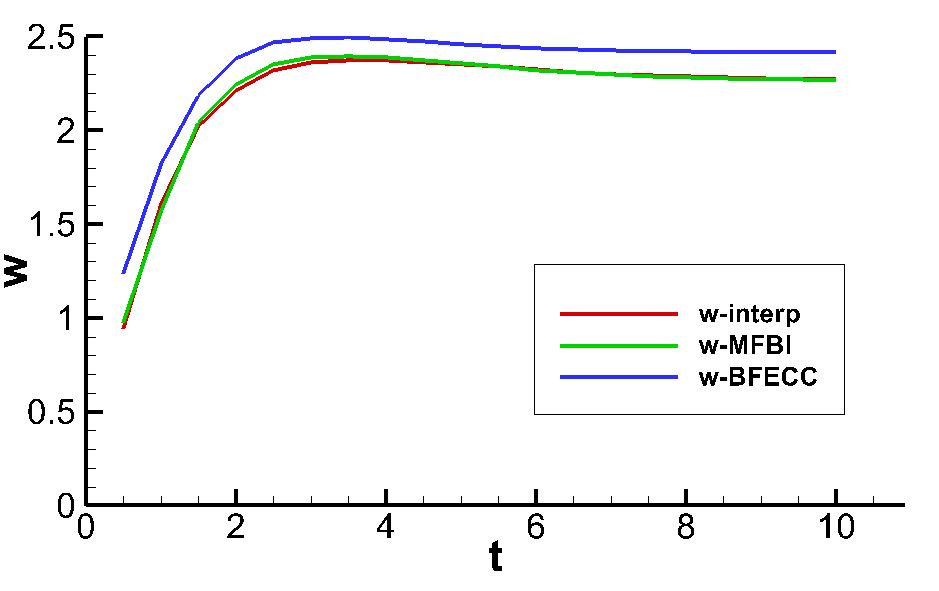}} \\  
  \caption{\small Accuracy order of corner flow at the interface, with inner cube parallel to the outer cube.}
  \label{corner_hor_conv_interf}
\end{figure}

\begin{figure}[!htbp]
 \centering
  \subfloat[p]{\label{hor_inner_ca_p}\includegraphics[width=0.35\textwidth]{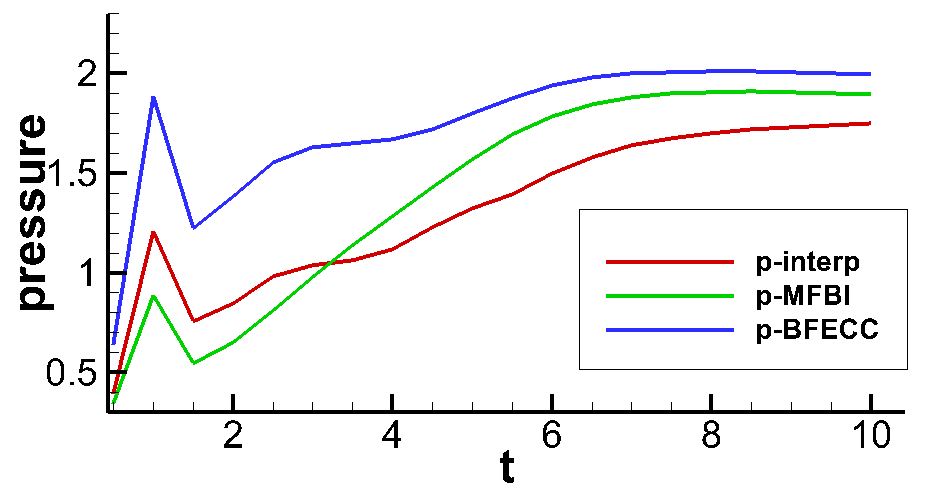}}  \ \
  \subfloat[u]{\label{hor_inner_ca_u}\includegraphics[width=0.35\textwidth]{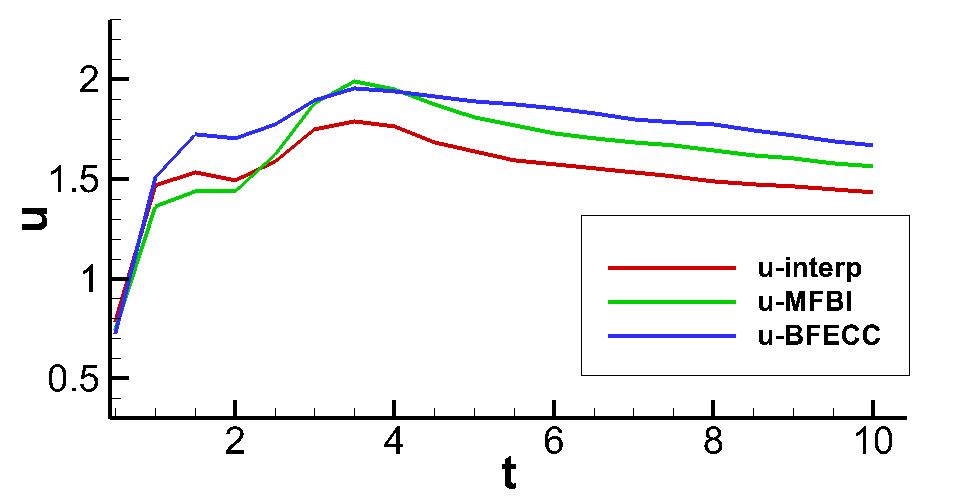}} \\
  \subfloat[v]{\label{hor_inner_ca_v}\includegraphics[width=0.35\textwidth]{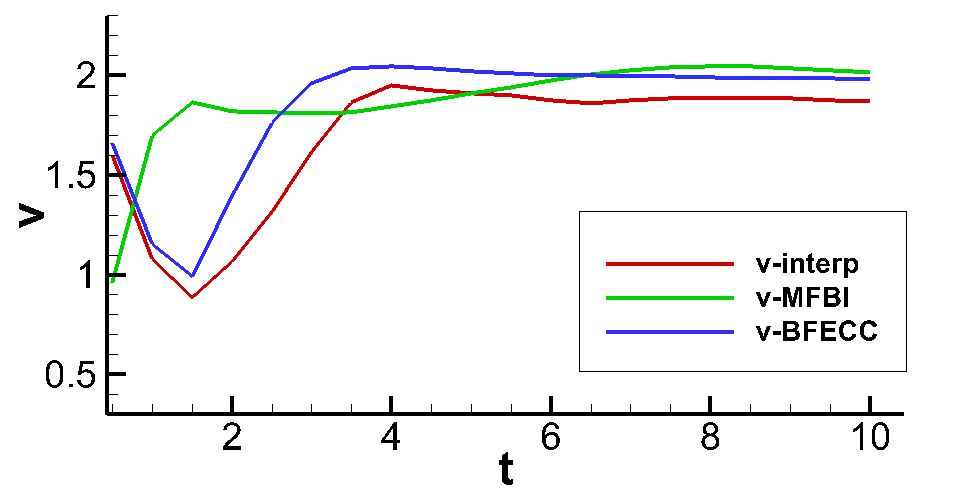}}  \ \
  \subfloat[w]{\label{hor_inner_ca_w}\includegraphics[width=0.35\textwidth]{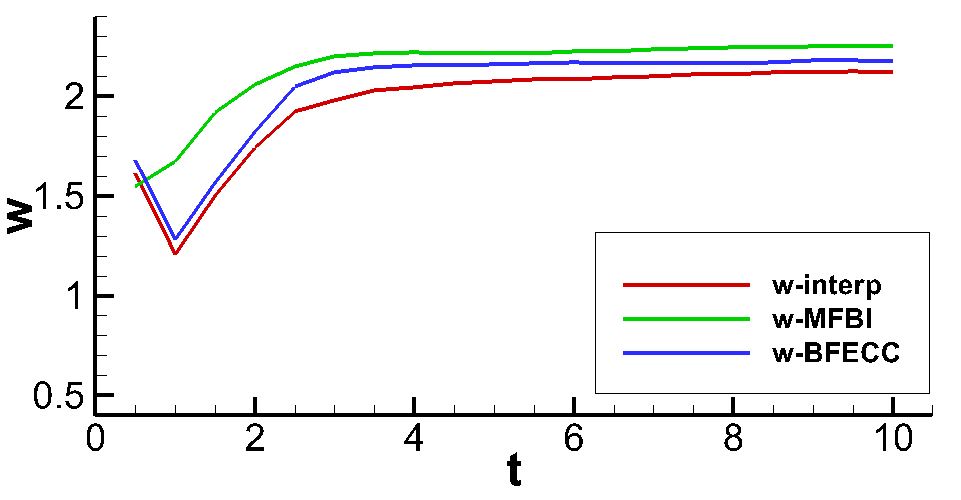}} \\  
  \caption{\small Accuracy order of corner flow in the inner cube, with inner cube parallel to the inner cube.}
  \label{corner_hor_conv_inner}
\end{figure}

\begin{figure}[!htbp]
 \centering
  \subfloat[p]{\label{hor_outer_ca_p}\includegraphics[width=0.35\textwidth]{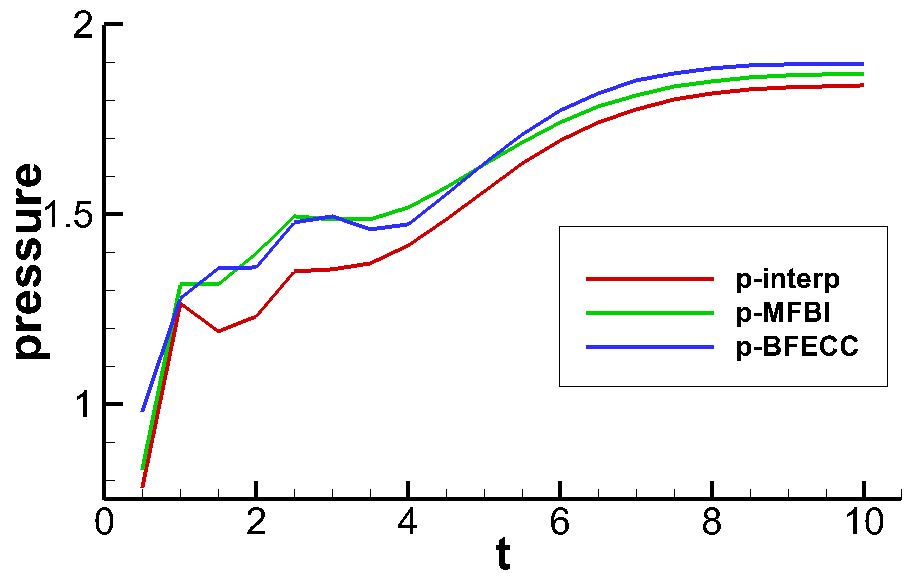}}  \ \
  \subfloat[u]{\label{hor_outer_ca_u}\includegraphics[width=0.35\textwidth]{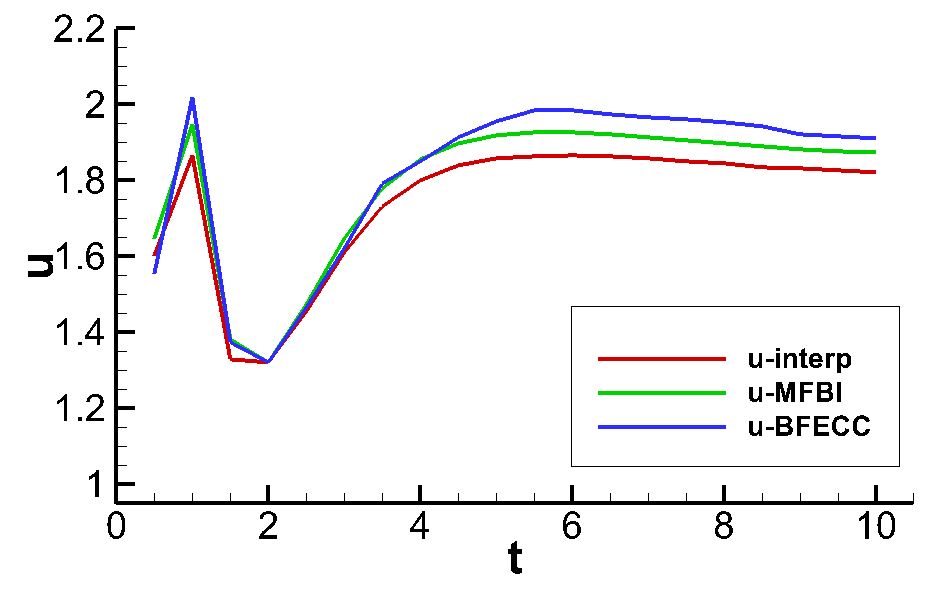}} \\
  \subfloat[v]{\label{hor_outer_ca_v}\includegraphics[width=0.35\textwidth]{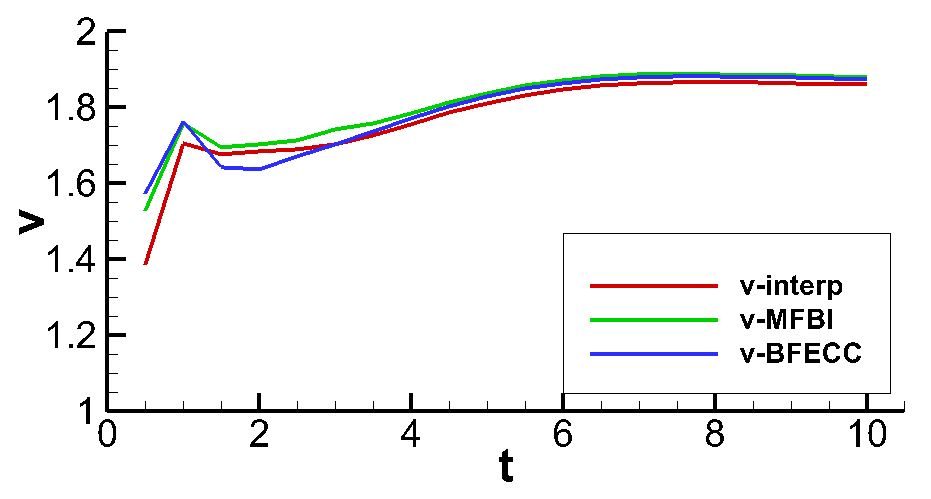}}  \ \
  \subfloat[w]{\label{hor_outer_ca_w}\includegraphics[width=0.35\textwidth]{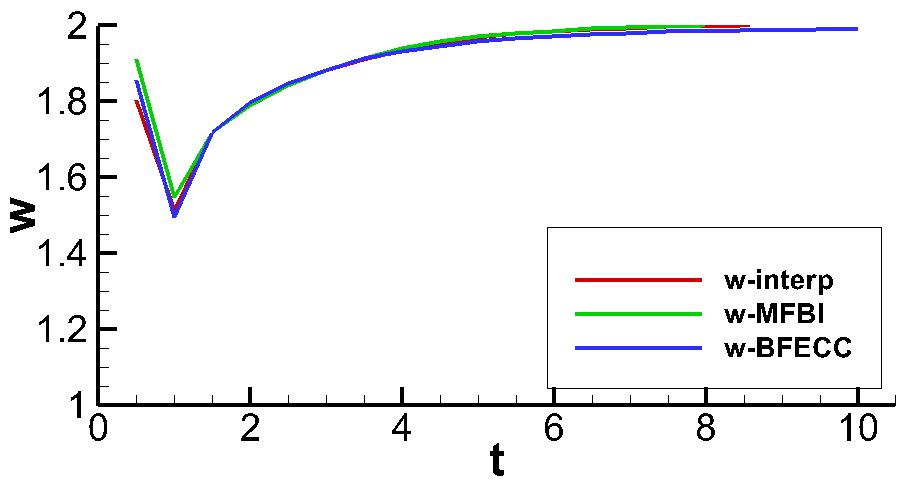}} \\  
  \caption{\small Accuracy order of corner flow in the outer cube,  with inner cube parallel to the outer cube.}
  \label{corner_hor_conv_outer}
\end{figure}

The estimated solution accuracy order, $\kappa $,  for the corner flow with the inner cube placed parallel to the outer one are plotted in Figs. \ref{corner_hor_conv_interf},  \ref{corner_hor_conv_inner}, and  \ref{corner_hor_conv_outer}, with show the accuracy order at the interface, in the inner cube, and in the outer cube, respectively. It is seen that the estimated accuracy orders are relatively low, and then after the flow develops for about 5-6 s, they become higher and stable. The low values at the initial time is attributed to that fact that the flow is highly transient and it is difficult for the formula (\ref {order_estimate}) to capture the accuracy order. The figures show that in all parts of the flow field, the accuracy order of BFECC is a little higher than that by the interpolation,  indicating that BFECC indeed enhances the accuracy of the interpolation and verifying the motion to apply BFECC at the interface. Interestingly, it is seen that MFBI also provides a better accuracy than the interpolation, primarily within the inner and outer cubes , although both of them are 1\ts{st} order accurate.

In order to get a better view on the accuracy of the three interface treatments, the average values of accuracy order over all time are summarized in Table.\ref {tab:corner_hor_c_ave_time}, which clearly indicates statistically the order for better accuracy is BFECC, MFBI, and linear interpolation.

\begin {table}[!htbp]
\begin{center}
\begin{tabular}{|l|l|l|l|l|l|l|}
\hline
      corner-hor        &           & p    & u    & v    & w    & puvw   \\ \hline
\multirow{3}{*}{Interp} & outer     & 1.54 & 1.74 & 1.78 & 1.91 & 1.74 \\ \cline{2-7} 
                        & inner     & 1.32 & 1.53 & 1.70 & 1.96 & 1.63 \\ \cline{2-7} 
                        & interface & 1.48 & 1.51 & 1.55 & 2.20 & 1.68 \\ \hline
\multirow{3}{*}{MFBI}   & outer     & 1.61 & 1.79 & 1.81 & 1.92 & 1.78 \\ \cline{2-7} 
                        & inner     & 1.41 & 1.63 & 1.88 & 2.14 & 1.77 \\ \cline{2-7} 
                        & interface & 1.49 & 1.52 & 1.58 & 2.20 & 1.70 \\ \hline
\multirow{3}{*}{BFECC}  & outer     & 1.63 & 1.82 & 1.79 & 1.91 & 1.79 \\ \cline{2-7} 
                        & inner     & 1.75 & 1.74 & 1.85 & 2.04 & 1.85 \\ \cline{2-7} 
                        & interface & 1.93 & 1.66 & 1.67 & 2.34 & 1.90 \\ \hline
\multicolumn{7}{|l|}{'corner-hor': inner block placed parallel to outer block}\\ \hline
\end{tabular}%
\caption {Component convergence rate averaged over time for corner-hor.}
\label{tab:corner_hor_c_ave_time}
\end{center}
\end{table}

It is noted that, as indicated in above section, BFECC has certain restrictions in relative rotation of grids. In order to verify such understanding, another case of corner flow simulated with rotated inner cube has also been conducted. The result shows that still BFECC increases the accuracy of interpolation, although not as much as in the horizontal situation.

\subsection{Cavity Flow} 

Cavity flow has been used to test numerical methods including interface algorithms \cite{Goda1979,Tang_2006_b}. Within the cavity, the fluid is initially stationary, and it flows due to a gradual motion of the lid in the horizontal direction:
\begin{equation}
\label{top_vel}
u=e^{-1/(5t)}+0.01. 
\end{equation}

\begin{figure}[!htbp]
\centering
\includegraphics[scale=0.3]{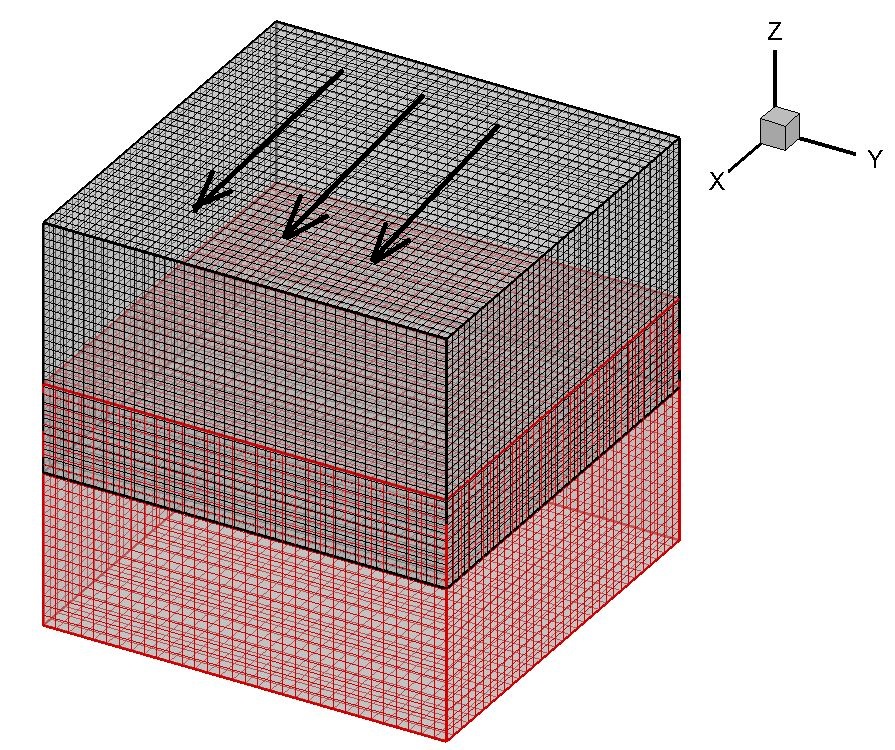}
\caption{Meshes for transient cavity flow}
\label{fig:cavity_mesh}
\end{figure}

\noindent As the flow is very complicated and it is difficult to estimate accuracy order of numerical solutions at high Reynolds numbers, we chose Re=1. In the computation, the flow domain is decomposed into two parts, upper part and lower part with an overlapping region in between, as shown in Fig.\ref{fig:cavity_mesh}. With origin of the coordinates at the center of the cavity, the upper part covers $(-0.5,0.5)\times (-0.5,0.5)\times (-0.12,0.5)$, and the lower part occupies $(-0.5,0.5)\times(-0.5,0.5)\times(-0.5,0.1)$. The interfaces of the two parts are at $z=-0.12$ and $z=0.1$. Again, three sets of grids are used to estimate the convergence rate: $9^3$ for the lower part and $11^3$ for the upper part as coarse grids, $17^3$ and $21^3$ as medium grids, and $33^3$ and $41^3$ as fine grids. Also, a reference solution is achieved by employing a fine mesh of $61^3$. 

\begin{figure}[!htbp]
 \centering
  \subfloat[pressure (y=0)]{\label{cavity_p}\includegraphics[width=0.35\textwidth]{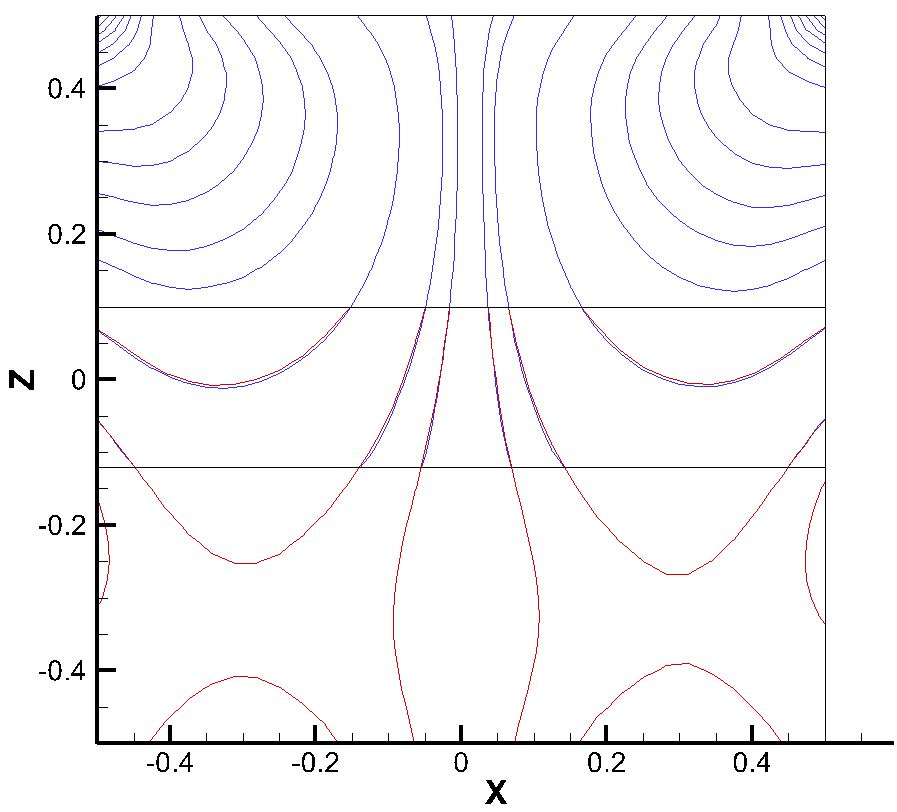}}  \ \
  \subfloat[u (y=0)]{\label{cavity_u}\includegraphics[width=0.35\textwidth]{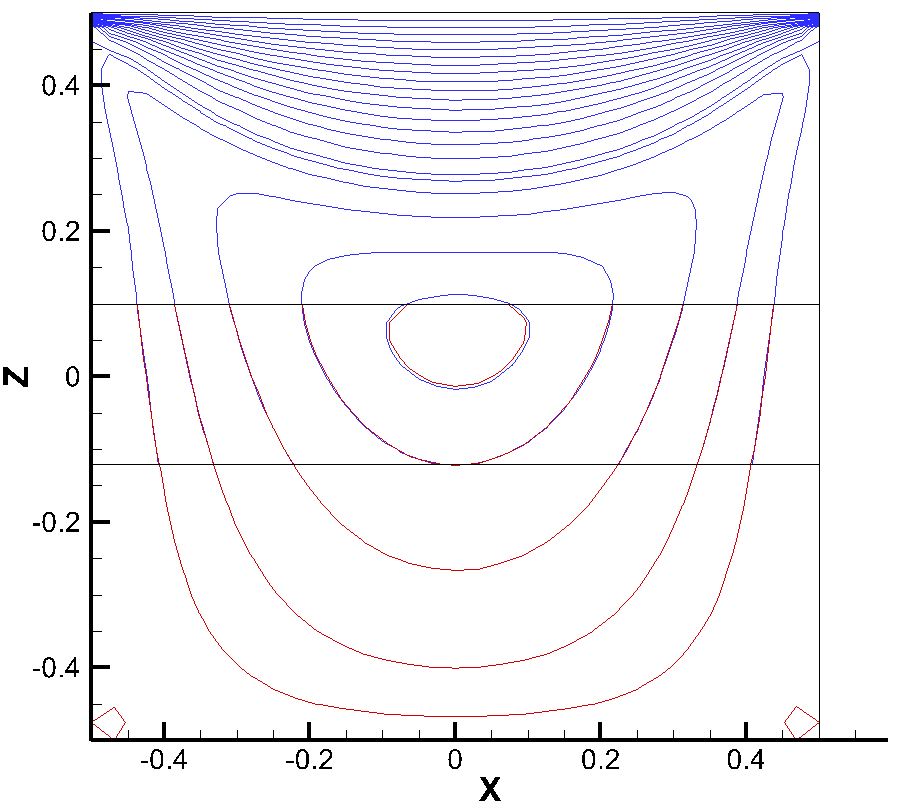}} \\
  \subfloat[v (x=0)]{\label{cavity_v}\includegraphics[width=0.35\textwidth]{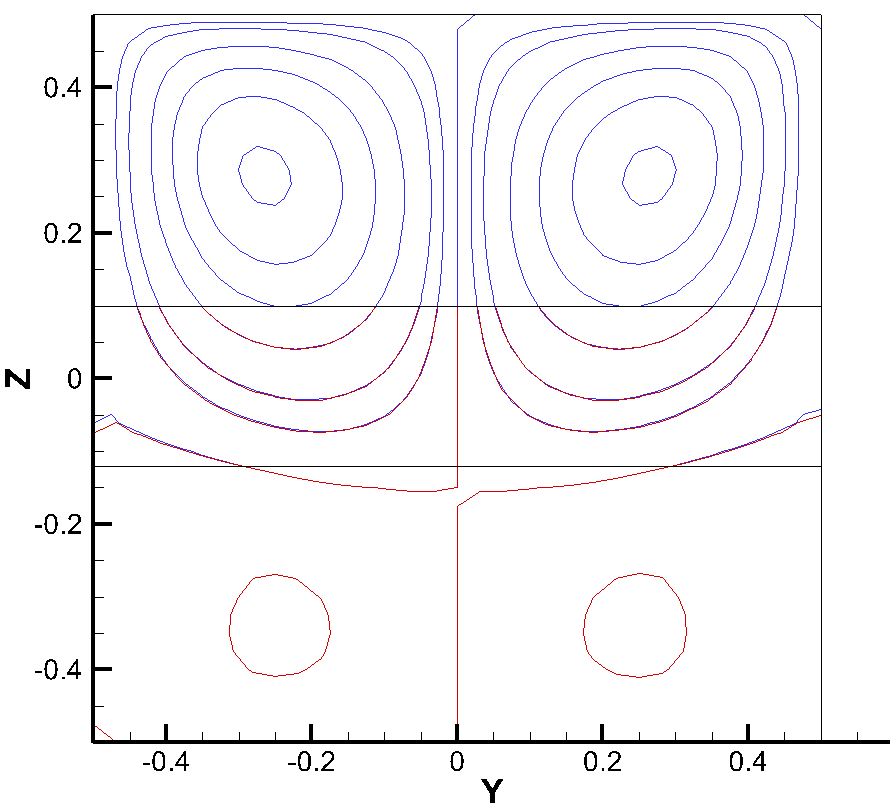}}  \ \
  \subfloat[w (y=0)]{\label{cavity_w}\includegraphics[width=0.35\textwidth]{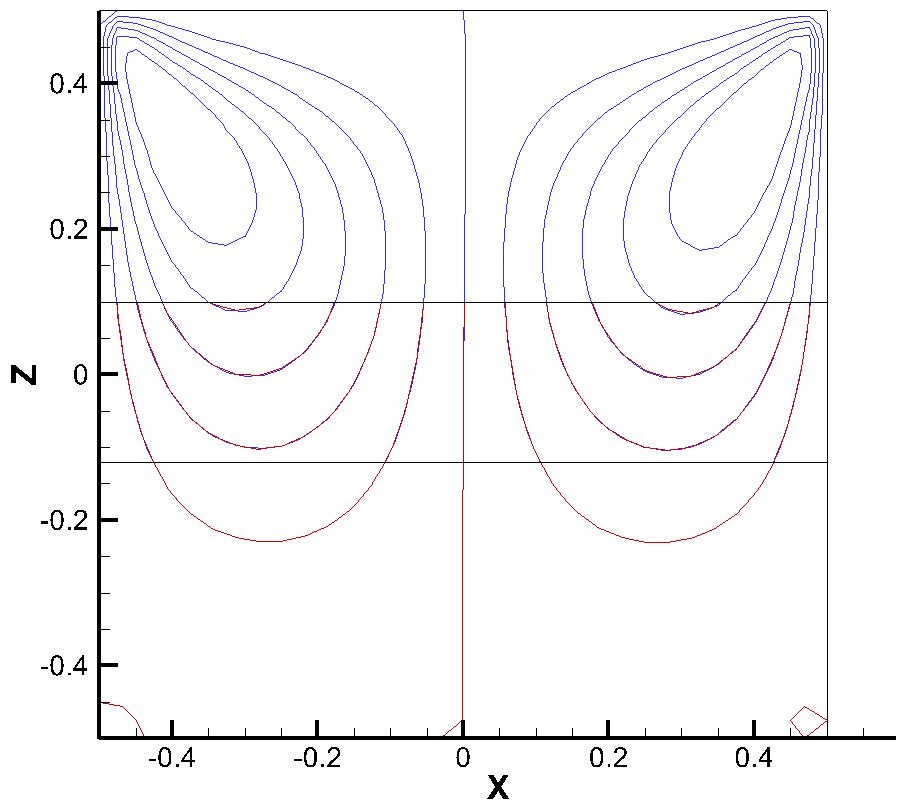}} \\  
  \caption{\small Solutions of cavity flow  obtained with BFECC, at t=1.0s}
  \label{cavity_flow_BFECC}
\end{figure}

\begin{figure}[!htbp]
 \centering
  \subfloat[p]{\label{cavity_ca_inter_p}\includegraphics[width=0.35\textwidth]{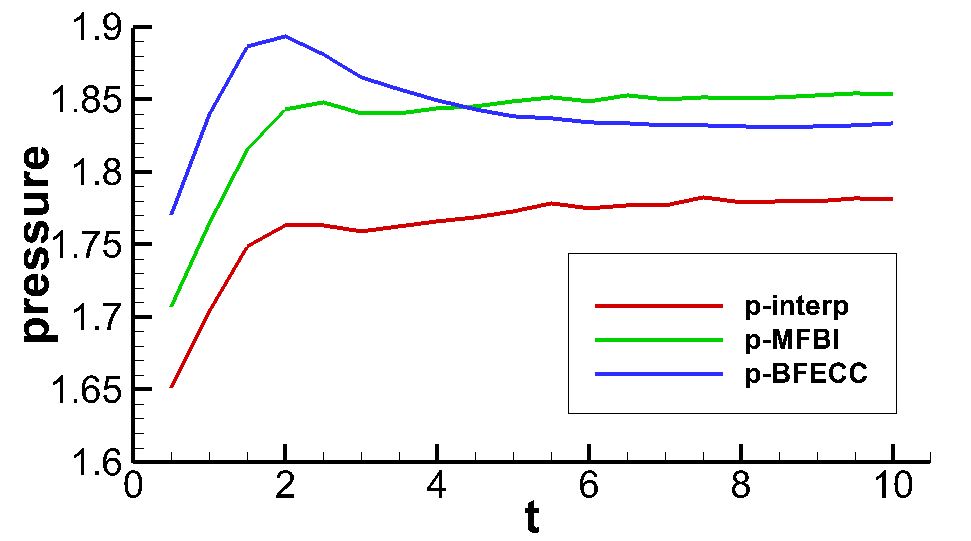}}  \ \
  \subfloat[u]{\label{cavity_ca_inter_u}\includegraphics[width=0.35\textwidth]{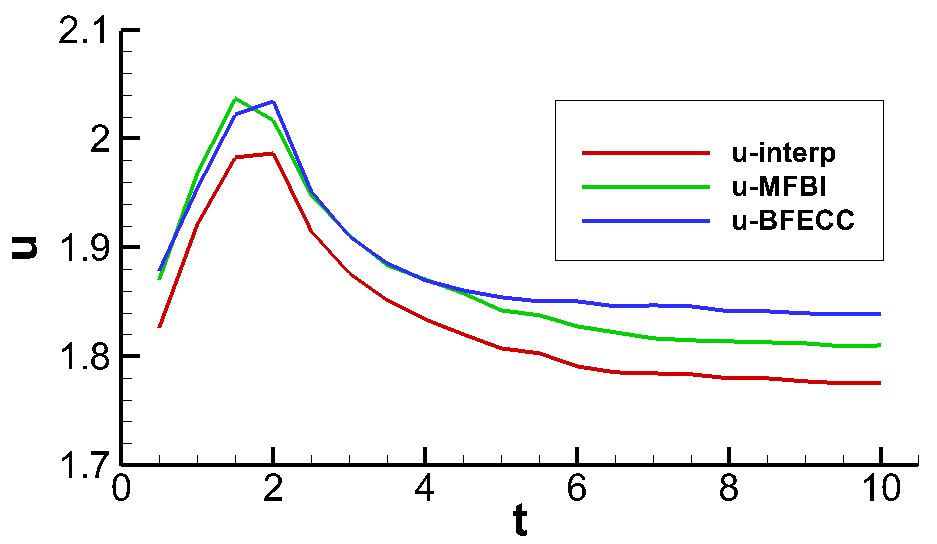}} \\
  \subfloat[w]{\label{cavity_ca_inter_v}\includegraphics[width=0.35\textwidth]{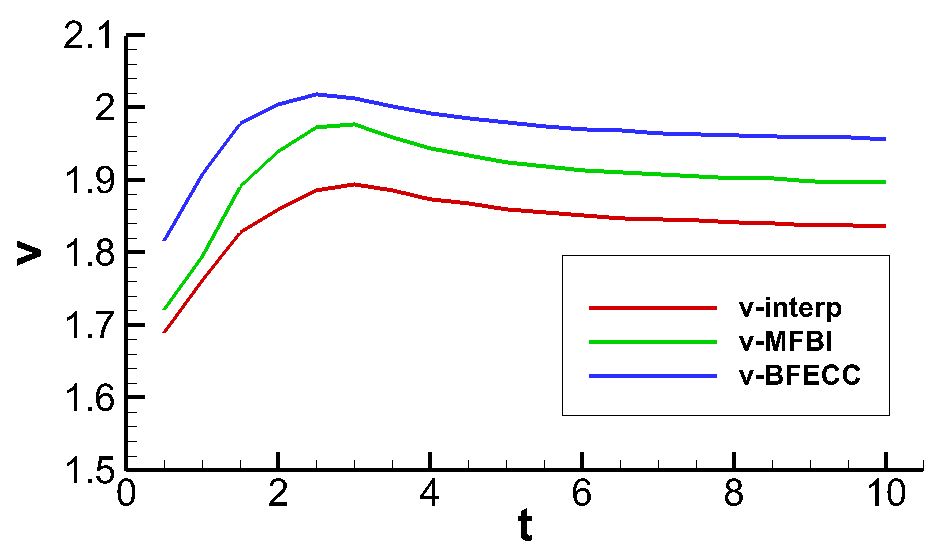}}  \ \
  \subfloat[v]{\label{cavity_ca_inter_w}\includegraphics[width=0.35\textwidth]{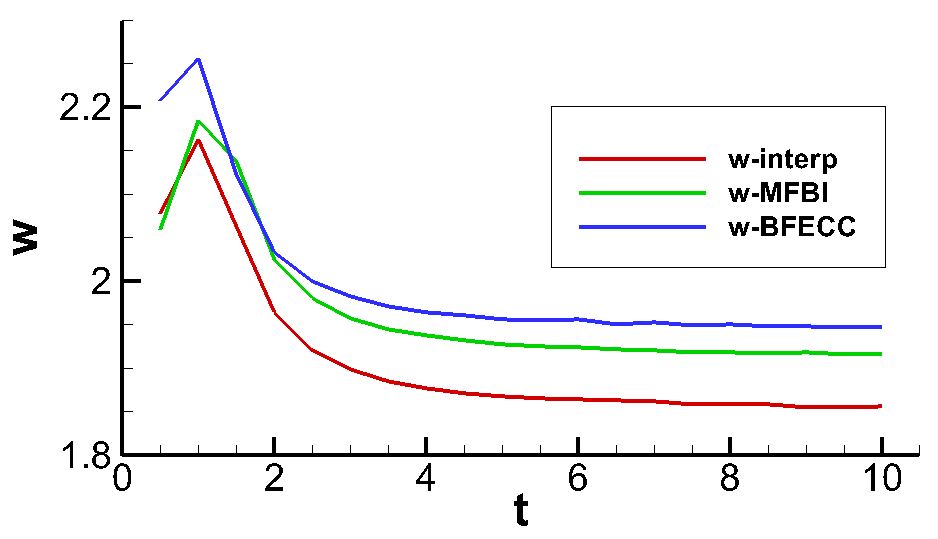}} \\  
  \caption{\small Accuracy order of cavity flow at the interface.}
  \label{cavity_conv_interf}
\end{figure}

\begin{figure}[!htbp]
 \centering
  \subfloat[p]{\label{cavity_ca_lower_p}\includegraphics[width=0.35\textwidth]{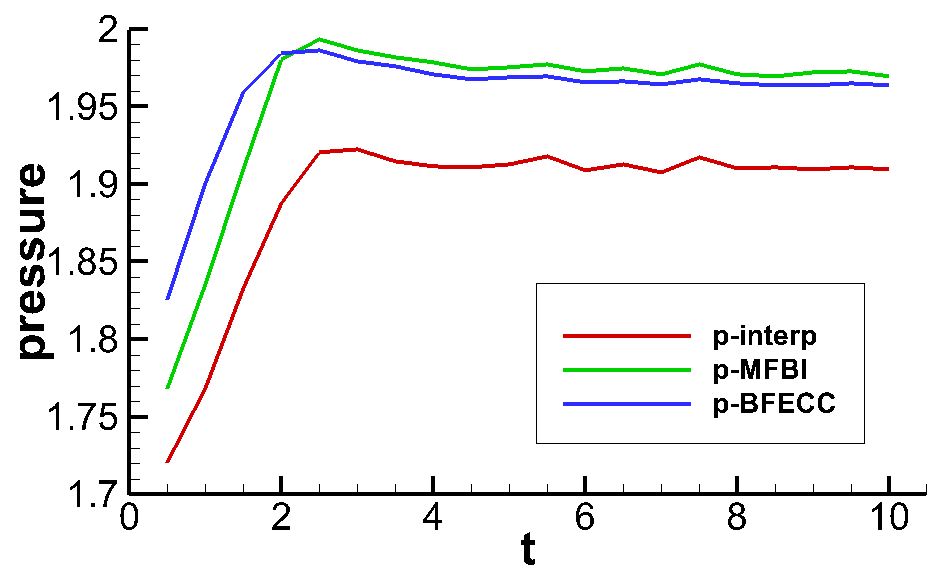}}  \ \
  \subfloat[u]{\label{cavity_ca_lower_u}\includegraphics[width=0.35\textwidth]{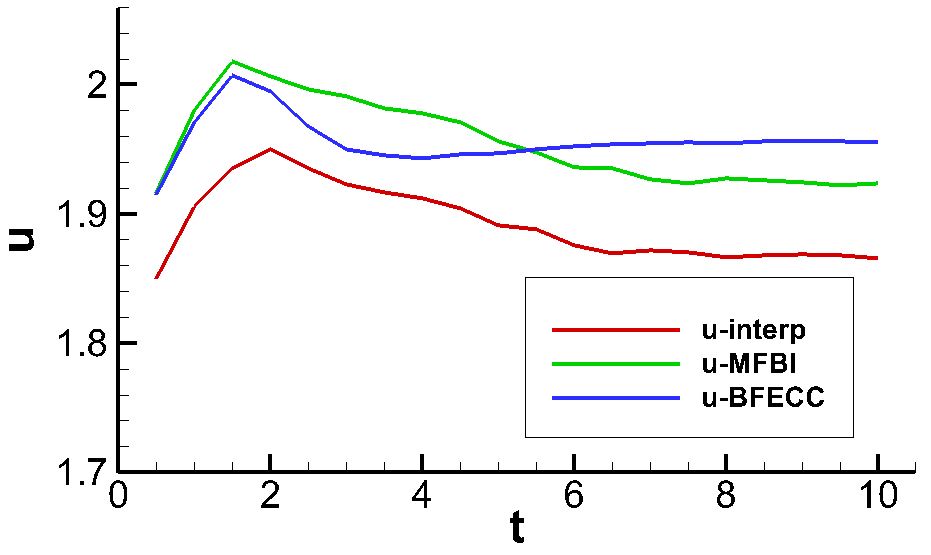}} \\
  \subfloat[w]{\label{cavity_ca_lower_v}\includegraphics[width=0.35\textwidth]{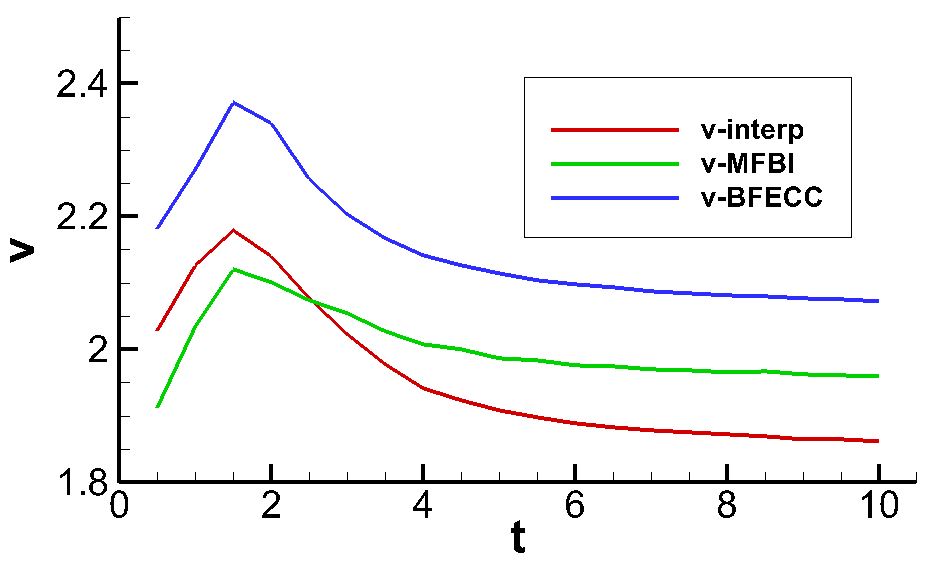}}  \ \
  \subfloat[v]{\label{cavity_ca_lower_w}\includegraphics[width=0.35\textwidth]{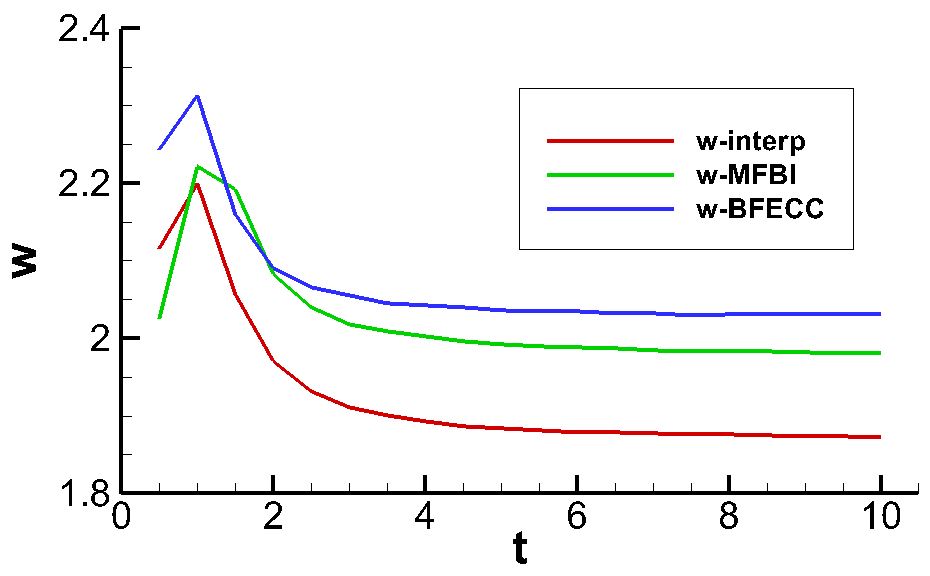}} \\  
  \caption{\small Accuracy order of cavity flow in the lower part.}
  \label{cavity_conv_lower}
\end{figure}

\begin{figure}[!htbp]
 \centering
  \subfloat[p]{\label{cavity_ca_up_p}\includegraphics[width=0.35\textwidth]{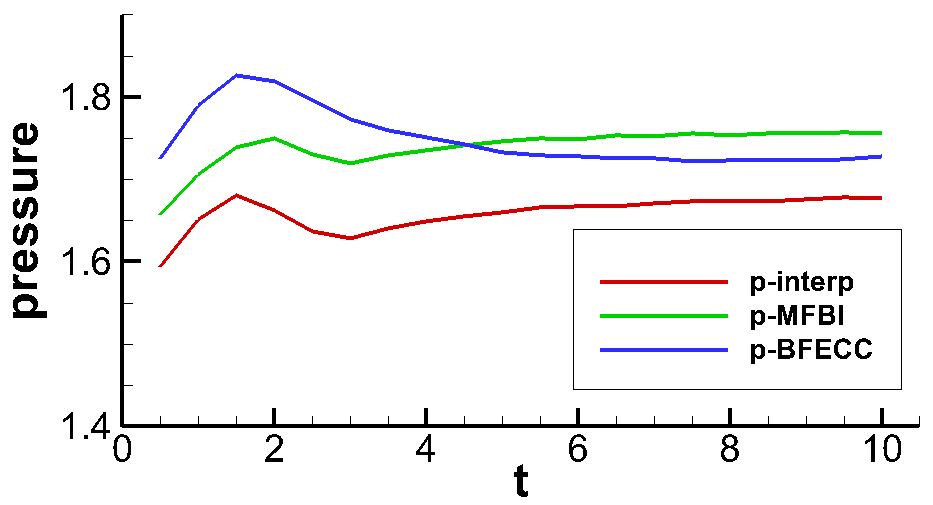}}  \ \
  \subfloat[u]{\label{cavity_ca_up_u}\includegraphics[width=0.35\textwidth]{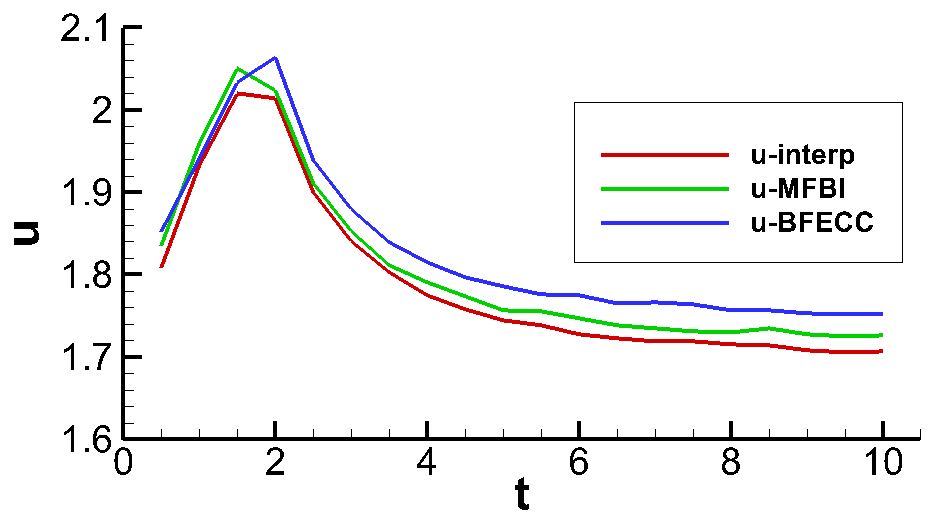}} \\
  \subfloat[w]{\label{cavity_ca_up_v}\includegraphics[width=0.35\textwidth]{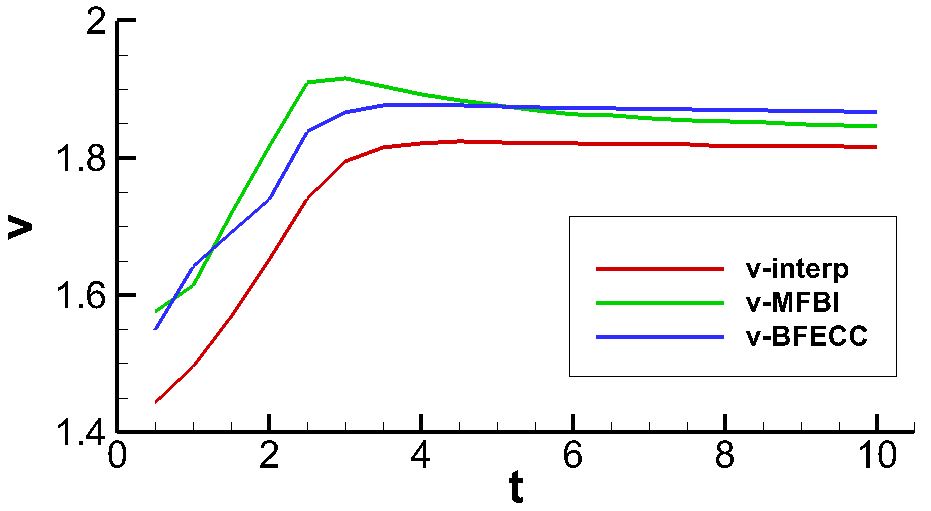}}  \ \
  \subfloat[v]{\label{cavity_ca_up_w}\includegraphics[width=0.35\textwidth]{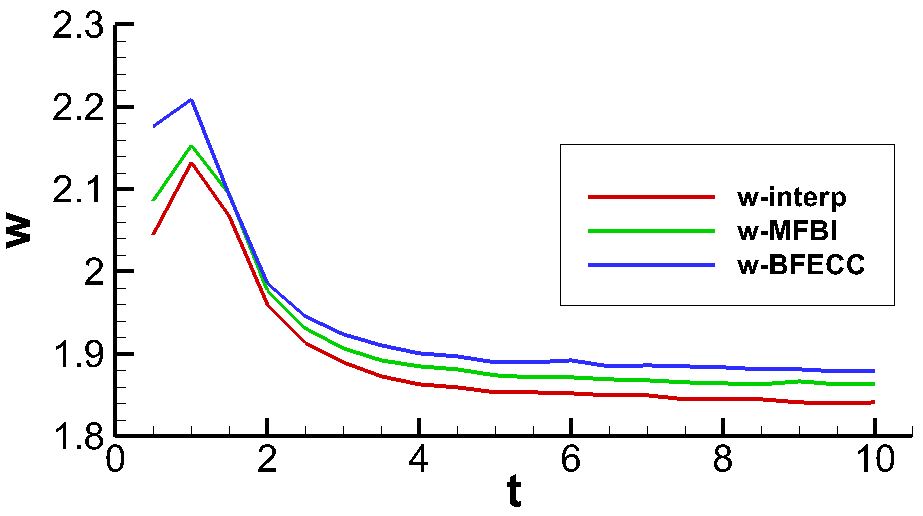}} \\  
  \caption{\small Accuracy order of cavity flow in the upper part.}
  \label{cavity_conv_upper}
\end{figure}

The simulated flow fields are shown in Fig.\ref {cavity_flow_BFECC}. As seen in the figure, the solutions in the lower and the upper parts are basically overlapping. To further analyze the results, similarly to discussion of the corner flow in above section, accuracy orders of the solutions are estimated and then plotted in Figs. \ref{cavity_conv_interf}, \ref{cavity_conv_lower}, and \ref{cavity_conv_upper}, which present the accuracy at the interface, in the lower part, and in the upper part. It is noted that the accuracy at the interface combines that at the two parts' interfaces, and the region $0.2<z<0.5$ is excluded when estimating accuracy since the solutions at the upper left and right corners are singular. It is seen in these figures that the estimated accuracy is low at the first a few seconds, and later it reaches a stable value around 2\ts{nd} order. The result indicates that BFECC indeed enhances the accuracy of the numerical solution. 

The average values of accuracy order over all time are summarized in Table.\ref {tab:cavity_c_ave_time}, which shows statistically the order for better accuracy is BFECC, MFBI, and linear interpolation.

\begin {table}[!htbp]
\begin{center}
\begin{tabular}{|l|l|l|l|l|l|l|}
\hline
cavity                  &           & p      & u      & v      & w      & puvw   \\ \hline
\multirow{3}{*}{Interp} & upper     & 1.66 & 1.79 & 1.76 & 1.90 & 1.78 \\ \cline{2-7} 
                        & lower     & 1.89 & 1.89 & 1.95 & 1.93 & 1.92 \\ \cline{2-7} 
                        & interface & 1.76 & 1.83 & 1.84 & 1.91 & 1.84 \\ \hline
\multirow{3}{*}{MFBI}   & upper     & 1.74 & 1.81 & 1.83 & 1.92 & 1.82 \\ \cline{2-7} 
                        & lower     & 1.96 & 1.95 & 2.00 & 2.02 & 1.98 \\ \cline{2-7} 
                        & interface & 1.84 & 1.87 & 1.91 & 1.96 & 1.89 \\ \hline
\multirow{3}{*}{BFECC}  & upper     & 1.75 & 1.83 & 1.83 & 1.94 & 1.84 \\ \cline{2-7} 
                        & lower     & 1.96 & 1.96 & 2.15 & 2.07 & 2.03 \\ \cline{2-7} 
                        & interface & 1.84 & 1.88 & 1.97 & 1.99 & 1.92 \\ \hline
\end{tabular}%
\caption {Component convergence rate averaged over time for cavity}
\label{tab:cavity_c_ave_time}
\end{center}
\end{table}


\newpage

\section{Conclusion}
\label{Sec: Conclusion}
BFECC interpolation is able to improve the order of accuracy
of the underlying linear interpolation method by one when 
interpolating from a rectangular grid to its translation, even if
a smooth perturbation is applied to the translation. If the translated grid points are located at centroids of the original grid,
then BFECC interpolation (based on linear interpolation) is able to
reach fourth order accuracy, a ``super-convergence'' phenomenon. These properties could be useful in high dimensional data and will be studied in the future. When interpolating from a rectangular grid to its rotation, BFECC interpolation will not improve the order of
the underlying linear interpolation, but will usually reduce the size of the linear interpolation error.
   
   
The investigation here is based on structured meshes. But the BFECC method also works on triangle meshes \cite{Hu_2014_b,Hu_2016_b}. In the next step, it is planned to apply this method to the coupled modeling system 'SIFOM-FVCOM' \cite {Tang_2014_b} (SIFOM is a solver for the Naver-Stokes equations, and FVCOM is a solver for the  geophysical fluid dynamics equations). But a new treatment is needed as the two coupled models use different types of meshes (structured mesh for SIFOM and unstructured mesh for FVCOM).

\vskip 0.1 true cm

\bibliographystyle{plain}
\bibliography{references}

\newpage
\begin{appendices}
All works in the appendices are conducted using 'Mathematica', version: 11.3.0.

\section{BFECC interpolation in 1D based on linear interpolation is 3\ts{rd} order accurate}
\label{1D_BFECC_proof}

This part is derived in a 1D situation. The grid spacing in each domain is assumed to be the same (dx1=dx2). 'g' is the source domain and 'f' is the target domain. The code is attached below:

\noindent 1. Taylor series:
\begin{lstlisting}
gin1= f-(1+a)*fx+(1+a)^2*fxx/2-(1+a)^3*fxxx/6+(1+a)^4*fxxxx/24
gi  = f-a*fx    +a^2*fxx/2    -a^3*fxxx/6    +a^4*fxxxx/24
gip1= f+(1-a)*fx+(1-a)^2*fxx/2+(1-a)^3*fxxx/6+(1-a)^4*fxxxx/24
gip2= f+(2-a)*fx+(2-a)^2*fxx/2+(2-a)^3*fxxx/6+(2-a)^4*fxxxx/24
\end{lstlisting}
\noindent 2. BFECC procedure:
\begin{lstlisting}
fin1= a*gi   + (1-a)*gin1
fi  = a*gip1 + (1-a)*gi
fip1= a*gip2 + (1-a)*gip1

ggi  = a*fin1 + (1-a)*fi
ggip1= a*fi   + (1-a)*fip1

gggi  = gi   + (gi-ggi)/2
gggip1= gip1 + (gip1-ggip1)/2

ffi   = a*gggip1 + (1-a)*gggi
\end{lstlisting}
\noindent The BFECC interpolation result:

$$ ffi = f + \frac{a(a-1)[(8a-4)fxxx-9a(a-1)fxxxx]}{24}$$

\noindent 3. The linear interpolation result:

$$ fi = f + \frac{a(1-a)[12fxx+(4-8a)fxxx+(1-3a+3a^2)fxxxx]}{24}$$
\newpage

\section{BFECC interpolation in 2D based on bilinear interpolation is 3\ts{rd} order accurate} 
\label{2D_BFECC_proof}

This part is derived in a 2D situation. The grid spacing in each domain is assumed to be the same (dx1=dx2, dy1=dy2, but dx1 $\neq$ dy1 \& dx2 $\neq$ dy2). 'g' is the source domain and 'f' is the target domain. The code is attached below:

\noindent 1. Taylor series:
\begin{lstlisting}
g1 =f -(1+b)*x*fx -(1+a)*y*fy +(1+b)^2*x^2*fxx/2 +(1+a)^2*y^2*fyy/2 
   +(1+a)*(1+b)*x*y*fxy -(1+b)^3*x^3*fxxx/6 -(1+a)^3*y^3*fyyy/6 
   -(1+b)^2*(1+a)*x^2*y*fxxy/2 -(1+b)*(1+a)^2*x*y^2*fxyy/2
g2 =f -b*x*fx     -(1+a)*y*fy +b^2*x^2*fxx/2     +(1+a)^2*y^2*fyy/2
   +(1+a)*b*x*y*fxy     -b^3*x^3*fxxx/6     -(1+a)^3*y^3*fyyy/6 
   -b^2*(1+a)*x^2*y*fxxy/2     -b*(1+a)^2*x*y^2*fxyy/2
g3 =f +(1-b)*x*fx -(1+a)*y*fy +(1-b)^2*x^2*fxx/2 +(1+a)^2*y^2*fyy/2 
   -(1+a)*(1-b)*x*y*fxy +(1-b)^3*x^3*fxxx/6 -(1+a)^3*y^3*fyyy/6 
   -(1-b)^2*(1+a)*x^2*y*fxxy/2 +(1-b)*(1+a)^2*x*y^2*fxyy/2
g4 =f +(2-b)*x*fx -(1+a)*y*fy +(2-b)^2*x^2*fxx/2 +(1+a)^2*y^2*fyy/2 
   -(1+a)*(2-b)*x*y*fxy +(2-b)^3*x^3*fxxx/6 -(1+a)^3*y^3*fyyy/6 
   -(2-b)^2*(1+a)*x^2*y*fxxy/2 +(2-b)*(1+a)^2*x*y^2*fxyy/2
g5 =f -(1+b)*x*fx -a*y*fy     +(1+b)^2*x^2*fxx/2 +a^2*y^2*fyy/2     
   +a*(1+b)*x*y*fxy     -(1+b)^3*x^3*fxxx/6 -a^3*y^3*fyyy/6     
   -(1+b)^2*a*x^2*y*fxxy/2     -(1+b)*a^2*x*y^2*fxyy/2
g6 =f -b*x*fx     -a*y*fy     +x^2*b^2*fxx/2     +y^2*a^2*fyy/2     
   +x*y*a*b*fxy         -x^3*b^3*fxxx/6     -y^3*a^3*fyyy/6     
   -x^2*y*b^2*a*fxxy/2         -x*y^2*b*a^2*fxyy/2
g7 =f +(1-b)*x*fx -a*y*fy     +(1-b)^2*x^2*fxx/2 +a^2*y^2*fyy/2     
   -(1-b)*a*x*y*fxy     +(1-b)^3*x^3*fxxx/6 -a^3*y^3*fyyy/6     
   -(1-b)^2*a*x^2*y*fxxy/2     +(1-b)*a^2*x*y^2*fxyy/2
g8 =f +(2-b)*x*fx -a*y*fy     +(2-b)^2*x^2*fxx/2 +a^2*y^2*fyy/2     
   -(2-b)*a*x*y*fxy     +(2-b)^3*x^3*fxxx/6 -a^3*y^3*fyyy/6     
   -(2-b)^2*a*x^2*y*fxxy/2     +(2-b)*a^2*x*y^2*fxyy/2
g9 =f -(1+b)*x*fx +(1-a)*y*fy +(1+b)^2*x^2*fxx/2 +(1-a)^2*y^2*fyy/2   
   -(1-a)*(1+b)*x*y*fxy -(1+b)^3*x^3*fxxx/6 +(1-a)^3*y^3*fyyy/6 
   +(1+b)^2*(1-a)*x^2*y*fxxy/2 -(1+b)*(1-a)^2*x*y^2*fxyy/2
g10=f -b*x*fx     +(1-a)*y*fy +b^2*x^2*fxx/2     +(1-a)^2*y^2*fyy/2 
   -b*(1-a)*x*y*fxy     -b^3*x^3*fxxx/6     +(1-a)^3*y^3*fyyy/6 
   +b^2*(1-a)*x^2*y*fxxy/2     -b*(1-a)^2*x*y^2*fxyy/2
g11=f +(1-b)*x*fx +(1-a)*y*fy +(1-b)^2*x^2*fxx/2 +(1-a)^2*y^2*fyy/2 
   +(1-a)*(1-b)*x*y*fxy +(1-b)^3*x^3*fxxx/6 +(1-a)^3*y^3*fyyy/6 
   +(1-b)^2*(1-a)*x^2*y*fxxy/2 +(1-b)*(1-a)^2*x*y^2*fxyy/2
g12=f +(2-b)*x*fx +(1-a)*y*fy +(2-b)^2*x^2*fxx/2 +(1-a)^2*y^2*fyy/2 
   +(1-a)*(2-b)*x*y*fxy +(2-b)^3*x^3*fxxx/6 +(1-a)^3*y^3*fyyy/6 
   +(2-b)^2*(1-a)*x^2*y*fxxy/2 +(2-b)*(1-a)^2*x*y^2*fxyy/2
g13=f -(1+b)*x*fx +(2-a)*y*fy +(1+b)^2*x^2*fxx/2 +(2-a)^2*y^2*fyy/2 
   -(2-a)*(1+b)*x*y*fxy -(1+b)^3*x^3*fxxx/6 +(2-a)^3*y^3*fyyy/6 
   +(1+b)^2*(2-a)*x^2*y*fxxy/2 -(1+b)*(2-a)^2*x*y^2*fxyy/2
g14=f -b*x*fx     +(2-a)*y*fy +b^2*x^2*fxx/2     +(2-a)^2*y^2*fyy/2 
   -(2-a)*b*x*y*fxy     -b^3*x^3*fxxx/6     +(2-a)^3*y^3*fyyy/6 
   +b^2*(2-a)*x^2*y*fxxy/2     -b*(2-a)^2*x*y^2*fxyy/2
g15=f +(1-b)*x*fx +(2-a)*y*fy +(1-b)^2*x^2*fxx/2 +(2-a)^2*y^2*fyy/2 
   +(2-a)*(1-b)*x*y*fxy +(1-b)^3*x^3*fxxx/6 +(2-a)^3*y^3*fyyy/6 
   +(1-b)^2*(2-a)*x^2*y*fxxy/2 +(1-b)*(2-a)^2*x*y^2*fxyy/2
g16=f +(2-b)*x*fx +(2-a)*y*fy +(2-b)^2*x^2*fxx/2 +(2-a)^2*y^2*fyy/2 
   +(2-a)*(2-b)*x*y*fxy +(2-b)^3*x^3*fxxx/6 +(2-a)^3*y^3*fyyy/6 
   +(2-b)^2*(2-a)*x^2*y*fxxy/2 +(2-b)*(2-a)^2*x*y^2*fxyy/2
\end{lstlisting}

\noindent 2. BFECC procedure:
\begin{lstlisting}
f1=(g1 *(1-a)+g5 *a)*(1-b)+(g2 *(1-a)+g6 *a)*b
f2=(g2 *(1-a)+g6 *a)*(1-b)+(g3 *(1-a)+g7 *a)*b
f3=(g3 *(1-a)+g7 *a)*(1-b)+(g4 *(1-a)+g8 *a)*b
f4=(g5 *(1-a)+g9 *a)*(1-b)+(g6 *(1-a)+g10*a)*b
f5=(g6 *(1-a)+g10*a)*(1-b)+(g7 *(1-a)+g11*a)*b
f6=(g7 *(1-a)+g11*a)*(1-b)+(g8 *(1-a)+g12*a)*b
f7=(g9 *(1-a)+g13*a)*(1-b)+(g10*(1-a)+g14*a)*b
f8=(g10*(1-a)+g14*a)*(1-b)+(g11*(1-a)+g15*a)*b
f9=(g11*(1-a)+g15*a)*(1-b)+(g12*(1-a)+g16*a)*b

gg6 =(f1*a+f4*(1-a))*b+(f2*a+f5*(1-a))*(1-b)
gg7 =(f2*a+f5*(1-a))*b+(f3*a+f6*(1-a))*(1-b)
gg10=(f4*a+f7*(1-a))*b+(f5*a+f8*(1-a))*(1-b)
gg11=(f5*a+f8*(1-a))*b+(f6*a+f9*(1-a))*(1-b)

ggg6 =g6 +0.5*(g6 -gg6 )
ggg7 =g7 +0.5*(g7 -gg7 )
ggg10=g10+0.5*(g10-gg10)
ggg11=g11+0.5*(g11-gg11)

ff5=(ggg6*(1-a)+ggg10*a)*(1-b)+(ggg7*(1-a)+ggg11*a)*b
\end{lstlisting}
\noindent The BFECC interpolation result:

$$ ffi = f + \frac{(b-3b^2+2b^3)dx^3*fxxx+(a-3a^2+2a^3)dy^3*fyyy}{6}$$

\noindent 3. The bilinear interpolation result:

$$ fi = f + \frac{(3b-3b^2)dx^2*fxx+(b-3b^2+2b^3)dx^3*fxxx}{6}$$

$$        + \frac{3a(1-a)dy^2*fyy+a(2a-1)(a-1)dy^3*fyyy}{6}$$

\begin{figure}[!htbp]
\centering
\includegraphics[scale=0.5]{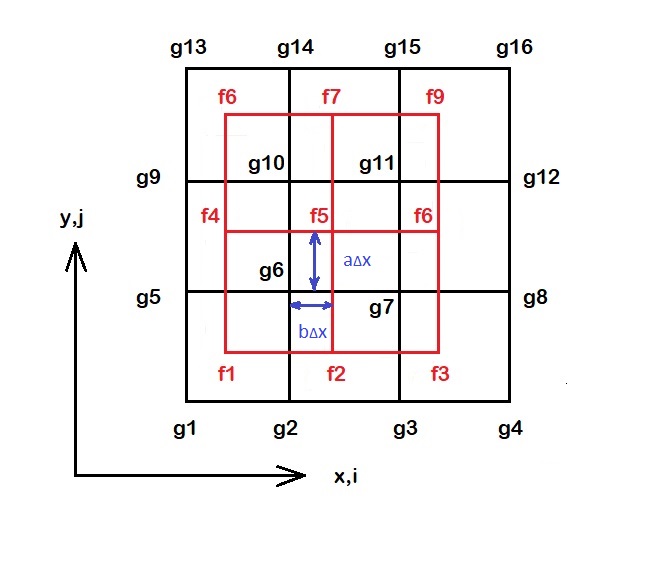}
\caption{Interpolation coefficients and node index}
\label{fig:2D_test1_domain}
\end{figure}

\newpage

\end{appendices}

\end{document}